\documentclass[12pt,reqno]{amsart}

\usepackage[utf8]{inputenc}
\usepackage[T1]{fontenc}
\usepackage[expansion=false]{microtype}

\usepackage{mathtools}          
\usepackage{amssymb,amsthm}
\usepackage{bm}                 

\usepackage[a4paper,margin=1in]{geometry}
\usepackage{enumitem}
\setlist[itemize]{leftmargin=*,itemsep=2pt,topsep=2pt}

\usepackage[dvipsnames,svgnames]{xcolor}
\usepackage[colorlinks=true,breaklinks=true,
            linkcolor=Mahogany,citecolor=ForestGreen,urlcolor=NavyBlue,
            pdfusetitle]{hyperref}

\theoremstyle{plain}
\newtheorem{theorem}{Theorem}[section]
\newtheorem{proposition}[theorem]{Proposition}
\newtheorem{lemma}[theorem]{Lemma}
\newtheorem{corollary}[theorem]{Corollary}
\newtheorem{theorema}{Theorem}

\numberwithin{equation}{section}
\allowdisplaybreaks

\newcommand{\R}{\mathbb{R}}
\newcommand{\bS}{\mathbb{S}}
\newcommand{\bD}{\mathrm{D}}
\newcommand{\bU}{\bm{U}}
\newcommand{\bV}{\bm{V}}
\newcommand{\bE}{\bm{E}}
\newcommand{\bW}{\bm{W}}
\DeclareMathOperator{\Tr}{Tr}
\DeclareMathOperator{\dist}{dist}
\DeclareMathOperator{\dv}{div}
\newcommand{\cH}{\mathcal{H}}
\newcommand{\cF}{\mathcal{F}}
\newcommand{\cT}{\mathcal{T}}
\newcommand{\cE}{\mathcal{E}}
\newcommand{\cG}{\mathcal{G}}

\newcommand{\C}{\mathbb{C}}
\newcommand{\bP}{\bm{P}}
\newcommand{\bQ}{\bm{Q}}
\newcommand{\bR}{\bm{R}}
\newcommand{\bZ}{\bm{Z}}
\newcommand{\bF}{\bm{F}}

\begin{document}

\title[Second order HUP for solenoidal fields]{Second order Heisenberg Uncertainty Principle for solenoidal fields: sharp constants, optimizers and sharp stability}

\author[A. X. Do]{Anh Xuan Do}
\address{Anh Xuan Do: Department of Mathematics, University of Connecticut, Storrs, CT
06269, USA}
\email{anh.do@uconn.edu}
\author[A. T. Duong]{Anh Tuan Duong}
\address{Anh Tuan Duong: Faculty of Mathematics and Informatics, Hanoi University of Science and Technology, Hanoi, Vietnam}
\email{tuan.duonganh@hust.edu.vn}
\author[N. Lam]{Nguyen Lam}
\address{Nguyen Lam: School of Science and the Environment, Grenfell Campus, Memorial
University of Newfoundland, Corner Brook, NL A2H 5G4, Canada}
\email{nlam@mun.ca}
\author[G. Lu]{Guozhen Lu}
\address{Guozhen Lu: Department of Mathematics, University of Connecticut, Storrs, CT
06269, USA}
\email{guozhen.lu@uconn.edu}
\author[V. H. Nguyen]{Van Hoang Nguyen}
\address{Van Hoang Nguyen: Department of Mathematics, FPT University, Ha Noi, Vietnam}
\email{vanhoang0610@yahoo.com,~hoangnv47@fe.edu.vn}
\date{\today}

\begin{abstract}
The second order Heisenberg Uncertainty Principle of Cazacu, Flynn and Lam is the uncertainty principle for the gradient of a function. Motivated by the question of Maz'ya [Integral Equations Operator Theory 90 (2018), Paper No. 25, Section 3.9] on the uncertainty principle for divergence-free vector fields, solved by Cazacu, Flynn and Lam for $N=2$ and by Hamamoto for $N\ge3$, we study the second order principle for divergence-free (solenoidal) vector fields $\bU$ on $\R^N$, $N\ge2$. We prove that
\[
\int_{\R^N}|\Delta\bU|^2dx\int_{\R^N}|x|^2|\nabla\bU|^2dx\ \ge\ C_N^2\Big(\int_{\R^N}|\nabla\bU|^2dx\Big)^2,
\]
where $C_N=1+\frac12\sqrt{N^2+8}$ for $N\ge3$ and $C_2=1+\sqrt2$. When $\bU$ is the velocity of an incompressible flow, $\int|\nabla\bU|^2$ and $\int|\Delta\bU|^2$ are, up to constant factors, the enstrophy and the palinstrophy of the flow. So the inequality bounds the palinstrophy of a flow from below by its enstrophy and by the spatial spread of its velocity gradient. The constant $C_N$ is sharp and larger than the constant $\frac{N+2}{2}$ for general vector fields, and all the optimizers are explicit fields with Kummer function profiles. We also obtain the sharp constants on the poloidal and the toroidal fields of each degree. Then we show that a divergence-free field which almost attains the equality is close to an optimizer, with sharp constants. In sum form, the deficit controls a chain of squared distances to an increasing sequence of explicit families of poloidal and toroidal fields, with all the constants sharp. In product form, we find the sharp stability constant in all dimensions. For $N\ge4$ it is $\frac12\big(\sqrt{N^2+4N-4}-\sqrt{N^2+8}\big)$, and it is attained by swirling flows, the toroidal fields of degree one. On the Fourier side, the sixth order problem for the poloidal fields of degree $k$ becomes a first order Caffarelli--Kohn--Nirenberg inequality for radial functions on $\R^{N+2k}$, and the problem for the toroidal fields becomes a weighted inequality on $\R^N$. Their stability follows from weighted Poincar\'e inequalities with Gaussian type measures.
\end{abstract}

\maketitle

\section{Introduction}\label{sec:intro}

\subsection{Stability of the Heisenberg Uncertainty Principle}
The main aim of this paper is the stability of the second order Heisenberg Uncertainty Principle (HUP) for solenoidal vector fields. The classical HUP states that
\begin{equation}\label{HUP}
\Big(\int_{\R^N}|\nabla u|^2dx\Big)\Big(\int_{\R^N}|x|^2|u|^2dx\Big)\ \ge\ \frac{N^2}{4}\Big(\int_{\R^N}|u|^2dx\Big)^2 ,
\end{equation}
and the constant $\frac{N^2}{4}$ is optimal and is attained by the Gaussian functions $\alpha e^{-\beta|x|^2}$, $\alpha\in\R$, $\beta>0$. We refer to \cite{FS97} for a survey on the uncertainty principle, and to \cite{Hei27, Wey31} for its origin in quantum mechanics. Uncertainty principles with fractional derivatives were studied in \cite{KPV22, Ste20}, and the HUP on half-spaces and orthants in \cite{LLLS25}. The stability question asks whether a function $u$ for which the two sides of \eqref{HUP} are almost equal must be close to the Gaussian functions, and how the distance to these functions is controlled by the deficit. We refer to \cite{Carlen} for a general discussion of stability for functional inequalities.

This question was first studied by McCurdy and Venkatraman \cite{MV21}. They proved that there are constants $C_1>0$ and $C_2(N)>0$ such that
\[
\int_{\R^N}|\nabla u|^2dx\int_{\R^N}|x|^2|u|^2dx-\frac{N^2}{4}\Big(\int_{\R^N}|u|^2dx\Big)^2\ \ge\ C_1\Big(\int_{\R^N}|u|^2dx\Big)d^2+C_2(N)\,d^4,
\]
where $d$ is the distance in $L^2(\R^N)$ from $u$ to the set of the Gaussian functions.
Their proof is by concentration-compactness, and it gives no information on $C_1$ and $C_2(N)$. Shortly after, Fathi \cite{Fat21} gave a short constructive proof with $C_1=\frac14$ and $C_2=\frac1{16}$. After the normalization $\int|u|^2=1$ and $\int|x|^2|u|^2=N$, the function $|u|^2$ is a probability density, \eqref{HUP} says that its Fisher information is at least $N$, and the difference is the relative Fisher information of $|u|^2$ with respect to the standard Gaussian measure. Fathi bounded this relative Fisher information from below by the Gaussian Poincar\'e inequality. His constants are explicit, but they are not sharp.

The sharp constants were found by Cazacu, Flynn, Lam and Lu in \cite{CFLL24}. Their approach has two steps. First, they use the deficit in product form
\[
\Big(\int_{\R^N}|\nabla u|^2dx\Big)^{\frac12}\Big(\int_{\R^N}|x|^2|u|^2dx\Big)^{\frac12}-\frac N2\int_{\R^N}|u|^2dx ,
\]
which is invariant under the dilations of $u$, and they prove that it is equal to
\[
\frac{\lambda^2}{2}\int_{\R^N}\Big|\nabla\Big(ue^{\frac{|x|^2}{2\lambda^2}}\Big)\Big|^2e^{-\frac{|x|^2}{\lambda^2}}dx,\qquad \lambda=\Big(\int_{\R^N}|x|^2|u|^2dx\Big/\int_{\R^N}|\nabla u|^2dx\Big)^{\frac14}.
\]
This identity shows at once that the optimizers of \eqref{HUP} are the Gaussian functions. Second, they apply to the function $ue^{|x|^2/(2\lambda^2)}$ the sharp Poincar\'e inequality for the Gaussian measure $e^{-|x|^2/\lambda^2}dx$, whose constant is $\frac{2}{\lambda^2}$. This gives the following result.

\begin{theorema}[\cite{CFLL24}]\label{ThA}
For all $u$ such that the integrals below are finite,
\[
\Big(\int_{\R^N}|\nabla u|^2dx\Big)^{\frac12}\Big(\int_{\R^N}|x|^2|u|^2dx\Big)^{\frac12}-\frac N2\int_{\R^N}|u|^2dx\ \ge\ \inf_{\alpha\in\R,\ \beta>0}\int_{\R^N}\big|u-\alpha e^{-\beta|x|^2}\big|^2dx .
\]
The inequality is sharp.
\end{theorema}

As a consequence, the sharp constants in the estimate of McCurdy and Venkatraman are $C_1=N$ and $C_2=1$ \cite{CFLL24}. Since the estimate of Theorem \ref{ThA} is attained by functions which are not Gaussian, one can ask again for its stability. Lam, Lu and Russanov \cite{LLR25} answered this question by combining the identity above with an improved Gaussian Poincar\'e inequality, in which the distance to the constants is replaced by the distance to the functions $c+d\cdot x$.

The approach of \cite{CFLL24} was extended by Do, Lam, Lu and Nguyen \cite{DLLN26} to the whole family of $L^2$-Caffarelli--Kohn--Nirenberg (CKN) inequalities \cite{CKN}
\begin{equation}\label{CKNfirst}
\Big(\int_{\R^N}\frac{|\nabla u|^2}{|x|^{2b}}dx\Big)^{\frac12}\Big(\int_{\R^N}\frac{|u|^2}{|x|^{2a}}dx\Big)^{\frac12}\ \ge\ C(N,a,b)\int_{\R^N}\frac{|u|^2}{|x|^{a+b+1}}dx ,
\end{equation}
which contains \eqref{HUP} as the case $a=-1$, $b=0$. The sharp constants and the optimizers of \eqref{CKNfirst} were found by Costa \cite{Cos08}, by Catrina and Costa \cite{CC09}, and by Cazacu, Flynn and Lam \cite{CFL21}, see also \cite{LL17}. Sharp CKN inequalities on stratified Lie groups and on Riemannian manifolds were obtained by Nguyen \cite{VHN18, VHN22}, and sharp weighted Sobolev and Gagliardo--Nirenberg inequalities were proved by mass transport in \cite{VHN15, VHN21}. In \cite{DLLN26}, the authors first consider the deficit in sum form, $\int\frac{|\nabla u|^2}{|x|^{2b}}+\int\frac{|u|^2}{|x|^{2a}}-2C(N,a,b)\int\frac{|u|^2}{|x|^{a+b+1}}$, which is not invariant under the dilations. For $1+b-a>0$ and $b\le\frac{N-2}{2}$, the sharp constant is $C(N,a,b)=\frac{N-1-a-b}{2}$, and the deficit is equal to
\[
\int_{\R^N}\Big|\nabla\Big(ue^{\frac{|x|^{1+b-a}}{1+b-a}}\Big)\Big|^2e^{-\frac{2|x|^{1+b-a}}{1+b-a}}\frac{dx}{|x|^{2b}} .
\]
So they need a Poincar\'e inequality for the measure $e^{-2|x|^{1+b-a}/(1+b-a)}|x|^{-2b}dx$. This measure is singular at the origin, and the classical proofs of the Gaussian Poincar\'e inequality, by the Ornstein--Uhlenbeck operator or by semigroup methods \cite{BGL14}, do not apply directly to it. They proved these weighted Poincar\'e inequalities, with their sharp constants and all their optimizers, by expansions in generalized Laguerre polynomials and spherical harmonics, and by a Kelvin type transform for the other ranges of $(a,b)$. The sharp constant is
\[
\min\Big\{2(1+b-a),\ \sqrt{(N-2-2b)^2+4(N-1)}-(N-2-2b)\Big\}.
\]
The first value comes from the radial functions, through the first Laguerre level, and the second one from the spherical harmonics of degree one. This gives the sharp stability of \eqref{CKNfirst} for all the parameters, together with the stability of the stability estimates. The deficit in product form is then treated by a dilation. See also \cite{DFLL23, LLLS26, LL26, LLR26} for related results.

\subsection{Solenoidal fields and second order inequalities}
The problem of finding the best constant in \eqref{HUP} when $u$ is replaced by a divergence-free vector field $\bU$ was raised by Maz'ya in \cite[Section 3.9]{Maz18}, with a motivation from hydrodynamics. More precisely, Maz'ya asked for the best constant $\mu_*(N)$ in the inequality
\[
\int_{\R^N}|\nabla\bU|^2dx\int_{\R^N}|x|^2|\bU|^2dx\ \ge\ \mu_*(N)\Big(\int_{\R^N}|\bU|^2dx\Big)^2,\qquad \bU\in C_0^\infty(\R^N)^N,\quad \dv\bU=0 .
\]
This question was answered by Cazacu, Flynn and Lam \cite{CFL22} in dimension two, and by Hamamoto \cite{Ham23} in all dimensions $N\ge3$:

\begin{theorema}[\cite{CFL22, Ham23}]\label{ThB}
We have $\mu_*(2)=4$ and $\mu_*(N)=\frac14\big(2+\sqrt{N^2-4N+12}\big)^2$ for $N\ge3$.
\end{theorema}

In dimension two, a divergence-free field is $\bU=\nabla^\perp\psi$ for a scalar function $\psi$, so that the problem becomes a second order inequality for $\psi$, which was solved in \cite{CFL22} by spherical harmonics. For $N\ge3$, Hamamoto \cite{Ham23} used the poloidal--toroidal decomposition of solenoidal fields. The toroidal part leads to first order problems, but the poloidal part leads to a one-dimensional minimization problem of the fourth order, and its solution is the main part of \cite{Ham23}. Simpler proofs of this one-dimensional problem were given in \cite{Ham24, HOT}, by writing the deficit as an integral of a square, in the spirit of the expanding-square method of \cite{Cos08, CFL23}. In \cite{DDLLNsol}, we gave a new proof of Theorem \ref{ThB} by a Fourier method, and we established the sharp stability of the HUP for solenoidal fields, together with chains of remainder terms. For instance, for $N\ge4$, the sharp stability constant of the HUP for solenoidal fields is $\frac12\big(N-\sqrt{N^2-4N+12}\big)$ \cite{DDLLNsol}. Solenoidal improvements of functional inequalities go back to Costin and Maz'ya \cite{CM08}, who improved the constant of the Hardy--Leray inequality \cite{Ler33} for axisymmetric solenoidal fields. We refer to \cite{Ham19, Ham20, Ham21, HT20, HT21} for the Hardy--Leray and Rellich--Leray inequalities for solenoidal and curl-free fields.

In \cite{CFL22}, Cazacu, Flynn and Lam also established the second order HUP
\begin{equation}\label{2ndHUPscalar}
\int_{\R^N}|\Delta u|^2dx\int_{\R^N}|x|^2|\nabla u|^2dx\ \ge\ \frac{(N+2)^2}{4}\Big(\int_{\R^N}|\nabla u|^2dx\Big)^2 ,
\end{equation}
and proved that the constant $\frac{(N+2)^2}{4}$ is sharp and is attained by the Gaussian functions. Writing $\bU=\nabla u$, \eqref{2ndHUPscalar} is the HUP for curl-free fields, see also \cite{CFL23}. Second order CKN inequalities were studied in \cite{CFL23, DN23}, and their stability in \cite{DN26}. The stability of \eqref{2ndHUPscalar} was first studied by Duong and Nguyen in \cite{DN25}. They wrote the deficit with the help of the quantity $\int\|\nabla^2v-x\otimes\nabla v\|^2e^{-|x|^2}dx$ and used the Hermite expansion of an Ornstein--Uhlenbeck type operator. Then Do, Lam and Lu \cite{DoLL26} used a different approach. They first replaced the product deficit by the deficit in sum form, which is easier to decompose. Then they applied the spherical harmonics decomposition, and they used an identity of \cite{DDLL25} which converts a radial second order problem on $\R^{N+2k}$ into a radial first order problem on $\R^{N+2k+2}$. At this point the strategy of \cite{CFLL24} applies, and they obtained the following result.

\begin{theorema}[\cite{DoLL26}]\label{ThC}
For all $u$ such that the integrals below are finite,
\begin{multline*}
\Big(\int_{\R^N}|\Delta u|^2dx\Big)^{\frac12}\Big(\int_{\R^N}|x|^2|\nabla u|^2dx\Big)^{\frac12}-\frac{N+2}{2}\int_{\R^N}|\nabla u|^2dx\\
\ge\ \frac{C(N)}{2}\inf_{\alpha\in\R,\ \beta>0}\int_{\R^N}\big|\nabla\big(u-\alpha e^{-\beta|x|^2}\big)\big|^2dx ,
\end{multline*}
and the constant $\frac{C(N)}{2}$ is sharp. Also, $\sqrt{N^2+4N-4}-N\le C(N)\le\frac{2N}{N+2}$.
\end{theorema}

In \cite{DLLNfrac}, Do, Lam, Lu and Nguyen showed that the constant $C(N)$ in Theorem \ref{ThC} is equal to $\sqrt{N^2+4N-4}-N$. Their observation is that, on the Fourier side, the second order HUP is a first order CKN inequality. Indeed, up to constant factors,
\[
\int_{\R^N}|\Delta u|^2dx=\int_{\R^N}|\xi|^4|\cF u|^2d\xi,\qquad \int_{\R^N}|\nabla u|^2dx=\int_{\R^N}|\xi|^2|\cF u|^2d\xi,\]
\[\int_{\R^N}|x|^2|\nabla u|^2dx=\int_{\R^N}|\xi|^2|\nabla\cF u|^2d\xi .
\]
So the weight $|x|^2$ acting on $\nabla u$ becomes again a first order quantity, and no derivative is lost. Hence the deficit of \eqref{2ndHUPscalar} for $u$ is the deficit of \eqref{CKNfirst} with $a=-2$ and $b=-1$ for $\cF u$, and the exact value of $C(N)$ follows from \cite{DLLN26}. In \cite{DLLNfrac}, the same idea also gives two families of sharp fractional HUPs, with their sharp constants, all their optimizers and their sharp stability, and chains of remainder terms whose constants are the successive spectral gaps of an explicit operator. See also \cite{HY25}.

So there are now two methods to obtain sharp stability constants for inequalities of HUP type. The first one, used in \cite{CFLL24, DLLN26, DoLL26, LLR25}, writes the deficit as an integral of a square, then factors out the optimizer and applies a Poincar\'e inequality for a Gaussian type measure. The second one, used in \cite{DLLNfrac, DDLLNsol}, applies the Fourier transform first and reduces the problem to a first order CKN inequality, to which the first method applies. In this paper we use the second method for the second order problem for solenoidal fields.

\subsection{Main results}
We now consider the second order HUP for vector fields. Applying \eqref{2ndHUPscalar} to each component of a vector field $\bU=(U_1,\dots,U_N)$ and using the Cauchy--Schwarz inequality, we see that
\begin{equation}\label{2ndHUP}
\int_{\R^N}|\Delta\bU|^2dx\int_{\R^N}|x|^2|\nabla\bU|^2dx\ \ge\ C^2\Big(\int_{\R^N}|\nabla\bU|^2dx\Big)^2
\end{equation}
holds for all vector fields with $C=\frac{N+2}{2}$. Here $|\nabla\bU|^2=\sum_i|\nabla U_i|^2$ and $|\Delta\bU|^2=\sum_i|\Delta U_i|^2$. The constant $\frac{N+2}{2}$ is sharp, and the optimizers are the fields whose components are proportional to a common Gaussian function. As in the problem of Maz'ya, we ask about the best constant in \eqref{2ndHUP} when $\bU$ is a solenoidal vector field, and about the stability of the corresponding inequality. The main purpose of this paper is to answer these two questions.

To state our results, we need some notation. For $k\ge0$, let $\cH_k$ be the space of spherical harmonics of degree $k$, and let $c_k:=k(k+N-2)$ be the corresponding eigenvalue of $-\Delta_{\bS^{N-1}}$. We identify $\phi\in\cH_k$ with the homogeneous harmonic polynomial $|x|^k\phi(x/|x|)$ of degree $k$ on $\R^N$, and we write $\cH_k^N$ for the space of vector fields whose components belong to $\cH_k$. For $N\ge2$ and $k\ge1$, we set
\begin{equation}\label{CNk}
C_P(N,k):=1+\frac12\sqrt{(N+2k-2)^2+8},\qquad C_T(N,k):=1+\frac12\sqrt{(N+2k-2)^2+4N-4}
\end{equation}
and
\begin{equation}\label{CN}
C_N:=\min\{C_P(N,1),\,C_T(N,1)\}.
\end{equation}
Since $\big(2C_T(N,k)-2\big)^2-\big(2C_P(N,k)-2\big)^2=4(N-3)$, we have $C_N=C_P(N,1)=1+\frac12\sqrt{N^2+8}$ for $N\ge3$, $C_2=C_T(2,1)=1+\sqrt2$, and $C_P(3,k)=C_T(3,k)$ for every $k\ge1$. We also set
\begin{equation}\label{betaD}
\beta_k:=\frac{\sqrt{(N+2k-2)^2+8}-(N+2k+2)}{2},\qquad D_k:=\sqrt{N^2+4c_k}-N ,
\end{equation}
so that $2C_P(N,k)=N+2k+4+2\beta_k$ and $2C_T(N,k)=N+2+D_k$. We use the Kummer functions
\begin{align}
p_k(r)&:={}_1F_1\Big(\frac{N+2k+\beta_k}{2},\ \frac{N+2k}{2},\ -\frac{r^2}{2}\Big),\label{gk}\\
t_k(r)&:={}_1F_1\Big(\frac{N+k}{2}+\frac{D_k}{4},\ \frac{N+2k}{2},\ -\frac{r^2}{2}\Big).\label{tk}
\end{align}
For $k\ge1$, we introduce the finite dimensional families
\begin{align*}
\cE_k&:=\Big\{\sum_{j=1}^k\bD(r^{j-1}p_j\phi_j):\phi_j\in\cH_j\Big\},\\
\cT_k&:=\Big\{\sum_{j=1}^kt_j(|x|)\bP_j(x):\bP_j\in\cH_j^N,\ x\cdot\bP_j(x)=0,\ \dv\bP_j=0\Big\}.
\end{align*}
Here $\bD$ is the poloidal generator of Hamamoto \cite{Ham23} recalled in Subsection \ref{sec:sol}, and in $\bD(r^{j-1}p_j\phi_j)$ the function $\phi_j$ is evaluated at $\sigma=x/|x|$ and $r=|x|$. For $\phi\in\cH_k$, the fields $\bD(r^{k-1}p_k\phi)$ are poloidal fields of degree $k$, and the fields $t_k(|x|)\bP(x)$ in the definition of $\cT_k$ are toroidal fields of degree $k$. For instance,
\[
\cT_1=\Big\{{}_1F_1\Big(\frac{D_1}{4}+\frac{N+1}{2},\ \frac{N+2}{2},\ -\frac{|x|^2}{2}\Big)Ax:\ A\in\R^{N\times N},\ A^T=-A\Big\}.
\]
For $N=2$, $\cT_1$ is spanned by $t_1(|x|)x^\perp$, where $x^\perp=(-x_2,x_1)$, and $\cT_k=\cT_1$ for every $k\ge1$. We use the dilations
\begin{equation}\label{dilation}
\bU_\lambda(x):=\lambda^{\frac N2-1}\bU(\lambda x),\qquad\lambda>0,
\end{equation}
which leave the right-hand side of \eqref{2ndHUP} invariant, and for a family $\cE$ we write $\cE^{\rm dil}:=\{\bE_\lambda:\bE\in\cE,\ \lambda>0\}$. All the fields above are smooth, and the integrals in \eqref{2ndHUP} are finite for them.

Our first main result gives the sharp constant of \eqref{2ndHUP} for solenoidal fields:

\begin{theorem}\label{thm:sharp}
Let $N\ge2$. Then for all solenoidal fields $\bU$,
\begin{equation}\label{2ndHUPsol}
\int_{\R^N}|\Delta\bU|^2dx\int_{\R^N}|x|^2|\nabla\bU|^2dx\ \ge\ C_N^2\Big(\int_{\R^N}|\nabla\bU|^2dx\Big)^2 ,
\end{equation}
where $C_N=1+\frac12\sqrt{N^2+8}$ for $N\ge3$ and $C_2=1+\sqrt2$. The constant $C_N$ is sharp. Moreover, $C_P(N,k)$ is the sharp constant on the poloidal fields of degree $k$, and, for $N\ge3$, $C_T(N,k)$ is the sharp constant on the toroidal fields of degree $k$. For $\bU\ne0$, the equality holds in \eqref{2ndHUPsol} if and only if $\bU\in\cE_1^{\rm dil}$ when $N\ge4$, $\bU\in(\cE_1+\cT_1)^{\rm dil}$ when $N=3$, and $\bU\in\cT_1^{\rm dil}$ when $N=2$.
\end{theorem}

The constant $C_N$ is strictly larger than the constant $\frac{N+2}{2}$ for general vector fields. For $N\ge4$ the optimizers are poloidal fields of degree one, and for $N=2$ they are toroidal fields. This is the same as for the first order problem \cite{CFL22, Ham23}. In dimension three, the poloidal and the toroidal constants are the same for each degree, and the optimizers are the sums of a poloidal optimizer and a toroidal optimizer at the same scale. In dimension two, every solenoidal field can be written as $\bU=\nabla^\perp\psi$, with $|\nabla\bU|=|\nabla^2\psi|$ and $|\Delta\bU|=|\nabla\Delta\psi|$. Hence \eqref{2ndHUPsol} with $N=2$ is a third order inequality for the scalar function $\psi$.

The inequality \eqref{2ndHUPsol} has a simple meaning for fluids. If $\bU$ is the velocity of an incompressible flow, its vorticity is the antisymmetric matrix $\omega=(\partial_iU_j-\partial_jU_i)_{i,j}$, which corresponds to $\mathrm{curl}\,\bU$ when $N=3$, and we write $|\omega|^2=\sum_{i,j}|\omega_{ij}|^2$. Since $\dv\bU=0$, an integration by parts gives
\[
\int_{\R^N}|\nabla\bU|^2dx=\frac12\int_{\R^N}|\omega|^2dx,\qquad \int_{\R^N}|\Delta\bU|^2dx=\frac12\int_{\R^N}|\nabla\omega|^2dx ,
\]
and these two quantities are, up to constant factors, the enstrophy of the flow \cite{MB02} and its palinstrophy. So \eqref{2ndHUPsol} bounds the palinstrophy of a flow from below by its enstrophy and by the spatial spread $\int|x|^2|\nabla\bU|^2dx$ of its velocity gradient.

Next, we study the stability of \eqref{2ndHUPsol}. We use the deficit
\begin{equation}\label{proddef}
\delta_{1,N}(\bU):=\Big(\int_{\R^N}|\Delta\bU|^2dx\Big)^{\frac12}\Big(\int_{\R^N}|x|^2|\nabla\bU|^2dx\Big)^{\frac12}-C_N\int_{\R^N}|\nabla\bU|^2dx ,
\end{equation}
which is invariant under the dilations \eqref{dilation}. As in \cite{DoLL26}, we will first study the stability with the following deficit:
\begin{equation}\label{sumdef}
\delta_{2,N}(\bU):=\int_{\R^N}|\Delta\bU|^2dx+\int_{\R^N}|x|^2|\nabla\bU|^2dx-2C_N\int_{\R^N}|\nabla\bU|^2dx .
\end{equation}
By the AM-GM inequality, $\delta_{2,N}(\bU)\ge2\delta_{1,N}(\bU)$, with equality when the first two integrals are equal. We measure the distances in the norm on the right-hand side of \eqref{2ndHUPsol}:
\[
\|\bU\|:=\Big(\int_{\R^N}|\nabla\bU|^2dx\Big)^{\frac12},\qquad
\dist(\bU,\cE):=\inf_{\bE\in\cE}\Big(\int_{\R^N}|\nabla(\bU-\bE)|^2dx\Big)^{\frac12} .
\]

Our second main result is a sharp stability estimate for \eqref{2ndHUPsol} with a chain of remainder terms:

\begin{theorem}\label{thm:sum}
For all solenoidal fields $\bU$, the following estimates hold.
\begin{itemize}
\item[(i)] Let $N\ge5$. Then
\begin{align*}
\delta_{2,N}(\bU)\ \ge\ &\big(\sqrt{N^2+4N-4}-\sqrt{N^2+8}\big)\dist^2(\bU,\cE_1)\\
&+\big(\sqrt{N^2+4N+12}-\sqrt{N^2+4N-4}\big)\dist^2(\bU,\cE_1+\cT_1)\\
&+\big(\sqrt{N^2+8N}-\sqrt{N^2+4N+12}\big)\dist^2(\bU,\cE_2+\cT_1)\\
&+\big(\sqrt{N^2+8N+24}-\sqrt{N^2+8N}\big)\dist^2(\bU,\cE_2+\cT_2)\\
&+\big(4+\sqrt{N^2+8}-\sqrt{N^2+8N+24}\big)\dist^2(\bU,\cE_3+\cT_2) .
\end{align*}
\item[(ii)] Let $N=4$. Then
\begin{align*}
\delta_{2,4}(\bU)\ \ge\ &\big(2\sqrt7-2\sqrt6\big)\dist^2(\bU,\cE_1)+\big(2\sqrt{11}-2\sqrt7\big)\dist^2(\bU,\cE_1+\cT_1)\\
&+\big(4\sqrt3-2\sqrt{11}\big)\dist^2(\bU,\cE_2+\cT_1)\\
&+\big(6\sqrt2-4\sqrt3\big)\dist^2(\bU,\cE_2+\cT_2)+\big(2\sqrt{19}-6\sqrt2\big)\dist^2(\bU,\cE_3+\cT_2)\\
&+\big(4+2\sqrt6-2\sqrt{19}\big)\dist^2(\bU,\cE_3+\cT_3) .
\end{align*}
\item[(iii)] Let $N=3$. Then
\begin{align*}
\delta_{2,3}(\bU)\ \ge\ &\big(\sqrt{33}-\sqrt{17}\big)\dist^2(\bU,\cE_1+\cT_1)+\big(\sqrt{57}-\sqrt{33}\big)\dist^2(\bU,\cE_2+\cT_2)\\
&+\big(4+\sqrt{17}-\sqrt{57}\big)\dist^2(\bU,\cE_3+\cT_3) .
\end{align*}
\item[(iv)] Let $N=2$. Then
\begin{align*}
\delta_{2,2}(\bU)\ \ge\ &\big(2\sqrt3-2\sqrt2\big)\dist^2(\bU,\cT_1)+\big(2\sqrt6-2\sqrt3\big)\dist^2(\bU,\cE_1+\cT_1)\\
&+\big(2\sqrt{11}-2\sqrt6\big)\dist^2(\bU,\cE_2+\cT_1)+\big(4+2\sqrt2-2\sqrt{11}\big)\dist^2(\bU,\cE_3+\cT_1) .
\end{align*}
\end{itemize}
In each case, all the constants are sharp, that is, none of them can be increased while the other ones are kept fixed. The equality holds when $\bU$ belongs to the largest family in the estimate.
\end{theorem}

The structure of these estimates is the following. The first family is the set of the optimizers of \eqref{2ndHUPsol} at a fixed scale. We decompose $\bU$ into its poloidal parts of degree $k$ and its toroidal parts of degree $k$. For each part, the number $2(C-C_N)$, where $C=C_P(N,k)$ or $C=C_T(N,k)$, is called the level of this part. The families become larger in the order of these levels. We add the extremal fields of the corresponding parts, and the constants are the differences between consecutive levels. The last level is $4$. It is the level of the first radial excitation of the optimizers (see Sections \ref{sec:pol} and \ref{sec:tor}). All the levels smaller than $4$ appear in the estimates. In particular, the sum of the constants in each estimate is $4$, and the first constant is the sharp stability constant for the deficit $\delta_{2,N}$.

Our third main result is the stability of \eqref{2ndHUPsol} with the deficit $\delta_{1,N}$:

\begin{theorem}\label{thm:prod}
Let $N\ge2$, let $\bU\ne0$ be a solenoidal field, and let
\[
\lambda_*:=\Big(\int_{\R^N}|x|^2|\nabla\bU|^2dx\Big/\int_{\R^N}|\Delta\bU|^2dx\Big)^{\frac14}.
\]
Then the following two estimates hold.
\begin{itemize}
\item[(i)] $\delta_{1,N}(\bU)$ is bounded from below by one half of the right-hand side of the estimate of Theorem \ref{thm:sum} for the field $\bU_{\lambda_*}$, and all the constants are sharp.
\item[(ii)] $\delta_{1,N}(\bU)$ is bounded from below by one half of the right-hand side of the estimate of Theorem \ref{thm:sum} for $\bU$, in which every family $\cE$ is replaced by $\cE^{\rm dil}$, and all the constants except the last one are sharp.
\end{itemize}
\end{theorem}

As a consequence, the sharp stability constant of \eqref{2ndHUPsol} with the deficit $\delta_{1,N}$ is
\[
\begin{cases}
\frac12\big(\sqrt{N^2+4N-4}-\sqrt{N^2+8}\big),&N\ge4,\\[2pt]
\frac12\big(\sqrt{33}-\sqrt{17}\big),&N=3,\\[2pt]
\sqrt3-\sqrt2,&N=2.
\end{cases}
\]
More precisely, $\delta_{1,N}(\bU)$ is bounded from below by this constant times the squared distance from $\bU$ to the set of all the optimizers of \eqref{2ndHUPsol}. For $N\ge4$, this constant comes from the toroidal fields of degree one, which are not optimizers, as for the first order problem \cite{DDLLNsol}. For $N=3$, it comes from the poloidal and the toroidal fields of degree two, which have the same level. For $N=2$, it comes from the poloidal fields of degree one. The formula for $N\ge4$ gives $0$ at $N=3$, since the toroidal fields of degree one become optimizers in dimension three.

We do not know if the last constant in Theorem \ref{thm:prod}(ii) is sharp. The test field which shows its sharpness in Theorem \ref{thm:prod}(i) is a radial excitation of an optimizer, and such a field is tangent to the set of the optimizers. The same happens for the scalar problem \cite{DLLNfrac}. The same argument, with the higher levels inside the blocks given by the Laguerre polynomials, gives further remainder terms above the level $4$. In this paper we stop at the first radial excitation, since all the remainder terms below this level are given by the explicit families $\cE_k$ and $\cT_k$.

\subsection{Our method}
Our approach follows the Fourier method that we developed in \cite{DDLLNsol} for the first order problem. The second order problem needs several new ingredients, and we describe them below, together with the difficulties that they solve. We think that the method can also be used for other inequalities of HUP type with differential constraints.

For general vector fields, the second order HUP \eqref{2ndHUP} follows from the scalar inequality \eqref{2ndHUPscalar} applied to each component. The divergence-free condition couples the components, and this simple argument is lost. The sharp constant becomes larger, and the optimizers are no longer Gaussian. As in \cite{Ham23, DDLLNsol}, we first use the poloidal--toroidal decomposition $\bU=\bU_P+\bU_T$ and the spherical harmonics decomposition. They remove the coupling, and the problem splits into a sum of independent problems, one for the poloidal part and one for the toroidal part of each degree. We call these parts the blocks.

Each block leads to a one-dimensional problem for a radial profile, and for the poloidal blocks this problem is of high order. Already for the first order HUP, Hamamoto obtained a minimization problem of the fourth order, and he wrote that the main difficulty of \cite{Ham23} lies in this problem. He solved it with an expansion in generalized Laguerre polynomials and by completing squares. In our case the integral $\int|\Delta\bU|^2$ of a poloidal field contains the third derivative of its radial profile, see \eqref{g3}, so that the problem is of the sixth order. Instead of solving it in the physical variable $r=|x|$, as in \cite{Ham23, Ham24, HOT}, we use Fourier analysis, as in \cite{DDLLNsol}. This is the main tool of our approach.

The Fourier transform lowers the order. Since $\int|x|^2|\nabla u|^2dx$ is equal to $\int|\xi|^2|\nabla\cF u|^2d\xi$ up to a constant factor, a weighted second order integral becomes a first order integral on the Fourier side, and no derivative is lost. This was observed in \cite{DLLNfrac} for scalar functions. To use it for a poloidal field $\bD(f_k\phi_k)$ of degree $k$, we write $f_k=r^{k-1}g_k$ and consider the radial function $G_k(X)=g_k(|X|)$ on $\R^{N+2k}$. Up to the same positive factor, the integrals $\int|\nabla\bU|^2$, $\int|\Delta\bU|^2$ and $\int|\nabla\bU|^2|x|^2$ are equal to
\[
\int|\Delta G_k|^2dX,\qquad \int|\nabla\Delta G_k|^2dX,\qquad \int|\Delta G_k|^2|X|^2dX+2\int|\nabla G_k|^2dX
\]
over $\R^{N+2k}$ (Proposition \ref{prop:lift}). So a poloidal field of degree $k$ in $\R^N$ behaves like a radial function in $\R^{N+2k}$, and we apply the Fourier transform on $\R^{N+2k}$. After the change of functions $\cF G_k=|\xi|^{\beta_k}F_k$, in which $\beta_k$ is chosen so that a zero order term cancels, the deficit of the block becomes
\[
\int|F_k|^2|\xi|^{6+2\beta_k}d\xi+\int|\nabla F_k|^2|\xi|^{4+2\beta_k}d\xi-2C_P(N,k)\int|F_k|^2|\xi|^{4+2\beta_k}d\xi ,
\]
up to a positive factor (see \eqref{blockw}). This is the deficit of a first order CKN inequality for $F_k$, and its optimizer is the Gaussian $e^{-|\xi|^2/2}$. The sixth order problem has become a first order problem, and the sharp constant $C_P(N,k)$ can be read off. The toroidal blocks are even simpler on the Fourier side. The Fourier transform exchanges the two conditions which define a toroidal field, since $\dv\bU=0$ becomes $\xi\cdot\cF\bU=0$ and $x\cdot\bU=0$ becomes $\dv\cF\bU=0$. So $\cF\bU_T$ is again toroidal, and with $\cF\bU_T=e^{-|\xi|^2/2}\bW_T$ we get
\[
(2\pi)^N\delta_{2,N}(\bU_T)=\int_{\R^N}|\nabla\bW_T|^2|\xi|^2e^{-|\xi|^2}d\xi-(2C_N-N-2)\int_{\R^N}|\bW_T|^2|\xi|^2e^{-|\xi|^2}d\xi
\]
(see \eqref{torw}). In each spherical degree, the right-hand side is the deficit of a weighted Poincar\'e inequality with a Gaussian type measure.

The Fourier side is also the right place to look for the optimizers. The approach of \cite{CFLL24, DoLL26} needs an identity in which the optimizer is factored out, and in the physical variable our optimizers have Kummer profiles, for which we do not have such an identity. On the Fourier side, the optimizers are Gaussian functions times a power of $|\xi|$, and the identity of \cite{CFLL24} can be used there. The inverse Fourier transform (Lemma \ref{lem:FT}) then gives all the extremal fields explicitly, in terms of the Kummer functions \eqref{gk} and \eqref{tk}. In this way we do not need to guess the extremals.

For the stability, the sharp constant is not enough. We need all the values of the deficit below some threshold, for the poloidal and the toroidal blocks of all degrees. Here the choice of the deficit is important. The natural deficit $\delta_{1,N}$ is invariant under the dilations, but it contains square roots of integrals, and it cannot be split over the blocks. As in \cite{CFLL24, DLLN26, DoLL26}, we first work with the deficit in sum form $\delta_{2,N}$, see \eqref{sumdef}, which is a quadratic form in $\bU$. This linearizes the problem: the stability in sum form becomes a spectral problem for the quadratic form $\int|\Delta\bU|^2+\int|\nabla\bU|^2|x|^2$ with respect to $\|\cdot\|^2$ on the solenoidal fields. On the Fourier side, the generalized Laguerre polynomials diagonalize this form and $\|\cdot\|^2$ in each block at the same time, which is not the case for the expansion of \cite{Ham23} in the physical variable. The values of $\delta_{2,N}$ relative to $\|\cdot\|^2$ are then the explicit numbers $2(C-C_N)+4j$, $j\ge0$, with $C=C_P(N,k)$ or $C=C_T(N,k)$, which we call levels. The bottom level $0$ gives the sharp constant and the optimizers, and the gaps between the levels give the stability constants. Theorem \ref{thm:sum} collects all the levels below $4$, each of them gives one term of the chain, and this is why all its constants are sharp. The level $4$ comes from the radial excitations of the optimizers. The order of the levels depends on the dimension, but it reduces to the comparison of the explicit numbers $c_k-N+3$ and $c_\ell$ (Lemma \ref{lem:levels}).

There is one more point for the toroidal blocks. The best approximation of $\bW_T$ in a given degree is first taken among all the polynomial fields of this degree, and these fields are not divergence-free in general. We show that this best approximation is again toroidal (Lemma \ref{lem:proj}). For degree one, this follows from a direct computation with the symmetric and the antisymmetric parts of a matrix. For higher degrees, we use a lemma on tangent divergence-free fields on the sphere (Lemma \ref{lem:sphere}). This step is needed to put the poloidal and the toroidal blocks together, since a field which is not toroidal is in general not orthogonal to the poloidal fields. It also gives remainder terms with solenoidal fields. We also use that the parts of degree zero of a toroidal field vanish. Without this fact, a field of degree zero would give the levels $N+2-2C_N+4j$, $j\ge0$. The level with $j=0$ is negative, and the level with $j=1$ is smaller than $4$, so the chains would stop earlier.

Finally, nothing is lost when we pass from $\delta_{1,N}$ to $\delta_{2,N}$. By the AM-GM inequality $\delta_{2,N}\ge2\delta_{1,N}$, and the two deficits are related exactly by
\[
\delta_{1,N}(\bU)=\frac12\,\delta_{2,N}(\bU_{\lambda_*}),
\]
where $\bU_{\lambda_*}$ is the dilate of $\bU$ for which the two integrals in the product are equal (see \eqref{transfer}). This gives Theorem \ref{thm:prod}, and the sharpness of the constants follows since the test fields of Theorem \ref{thm:sum} are balanced (Lemma \ref{lem:balanced}). Apart from the weighted Poincar\'e inequalities of \cite{DLLN26}, the proof uses only identities and elementary inequalities, and we never use a contradiction argument.

Our paper is organized as follows. In Section \ref{sec:pre}, we recall the weighted Poincar\'e inequalities with Gaussian type measures, some facts about the Fourier transform and the poloidal--toroidal decomposition. The poloidal fields are studied in Section \ref{sec:pol}, and the toroidal fields in Section \ref{sec:tor}. Theorem \ref{thm:sum} and the inequality in Theorem \ref{thm:sharp} are proved in Section \ref{sec:proofs}. Finally, Theorem \ref{thm:prod} and the equality cases in Theorem \ref{thm:sharp} are proved in Section \ref{sec:prod}.

Throughout this paper, the solenoidal fields are smooth and the three integrals in \eqref{2ndHUPsol} are finite, or they belong to the completion of such fields for the corresponding norm. The identities below are first proved for smooth fields with compact support, and they extend to this space by density. We write $r=|x|$ and $\sigma=x/r\in\bS^{N-1}$, and $\rho=|\xi|$ for the dual variable of the Fourier transform. For a radial function we use the same letter for the function and its profile. The fields are real valued, but we work with complex valued functions on the Fourier side. All the one-dimensional inequalities below hold for complex valued functions, by applying them to the real and imaginary parts.

\section{Preliminaries}\label{sec:pre}

\subsection{Weighted Poincar\'e inequalities with Gaussian type measures}
We use the weighted Poincar\'e inequalities with Gaussian type measures of Do, Lam, Lu and Nguyen \cite{DLLN26}. The first one is \cite[Theorem 2.1]{DLLN26} after the change of variables $s=\sqrt2r$.

\begin{theorema}[{\cite{DLLN26}}]\label{TA}
Let $\alpha>0$. Then for all $v\in W^{1,2}\big((0,\infty),r^{\alpha-1}e^{-r^2}dr\big)$,
\begin{equation}\label{GP}
\int_0^\infty|v'|^2r^{\alpha-1}e^{-r^2}dr\ \ge\ 4\inf_{c\in\C}\int_0^\infty|v-c|^2r^{\alpha-1}e^{-r^2}dr .
\end{equation}
The constant $4$ is sharp, and equality holds if and only if $v=a_0+a_1r^2$. More precisely, the functions $L_j^{(\frac\alpha2-1)}(r^2)$, $j\ge0$, where $L_j^{(\cdot)}$ are the generalized Laguerre polynomials \cite{Szego}, form an orthogonal basis of $L^2(r^{\alpha-1}e^{-r^2}dr)$, and the level $j$ contributes $4j$ times its squared norm to the left-hand side of \eqref{GP}. In particular the second level is spanned by $\frac\alpha2-r^2$.
\end{theorema}

We also use a Gaussian gauge. Let $m\ge1$, $p>-m$, and let $V=e^{-|\xi|^2/2}W$ on $\R^m$. Then $\nabla V=e^{-|\xi|^2/2}(\nabla W-\xi W)$, and an integration by parts of the cross term, using $\dv(\xi e^{-|\xi|^2}|\xi|^p)=(m+p-2|\xi|^2)e^{-|\xi|^2}|\xi|^p$, gives
\begin{equation}\label{gauss}
\int_{\R^m}|\nabla V|^2|\xi|^p\,d\xi+\int_{\R^m}|V|^2|\xi|^{p+2}\,d\xi-(m+p)\int_{\R^m}|V|^2|\xi|^p\,d\xi=\int_{\R^m}|\nabla W|^2e^{-|\xi|^2}|\xi|^p\,d\xi .
\end{equation}
This identity holds for vector-valued functions as well, componentwise. Combining \eqref{gauss} with Theorem \ref{TA}, we obtain the following estimate, which contains the radial case of \cite[Theorem 3.2]{DLLN26}.

\begin{corollary}\label{cor:CKN}
Let $m\ge1$, $p>-m$, and let $V$ be a radial function on $\R^m$ such that the integrals below are finite. Then
\begin{multline*}
\int_{\R^m}|\nabla V|^2|\xi|^p\,d\xi+\int_{\R^m}|V|^2|\xi|^{p+2}\,d\xi-(m+p)\int_{\R^m}|V|^2|\xi|^p\,d\xi\\
\ge\ 4\inf_{c\in\C}\int_{\R^m}\big|V-ce^{-\frac{|\xi|^2}{2}}\big|^2|\xi|^p\,d\xi .
\end{multline*}
The constant $4$ is sharp, and equality holds if and only if $V=(a_0+a_1|\xi|^2)e^{-|\xi|^2/2}$.
\end{corollary}

\begin{proof}
Let $W=e^{|\xi|^2/2}V$. By \eqref{gauss}, the left-hand side is $|\bS^{m-1}|\int_0^\infty|W'|^2\rho^{m+p-1}e^{-\rho^2}d\rho$. We apply Theorem \ref{TA} with $\alpha=m+p>0$, and we use that $|W-c|^2e^{-|\xi|^2}=|V-ce^{-|\xi|^2/2}|^2$.
\end{proof}

The next lemma is \cite[(2.7) and (2.8)]{DLLN26} with $a=-2$ and $b=-1$. We recall the short proof, since we also need the second level.

\begin{lemma}\label{lem:degree}
Let $N\ge2$, $\ell\ge1$, $Y\in\cH_\ell$, and let $W(\xi)=w(|\xi|)Y(\xi/|\xi|)$ on $\R^N$ be such that the left-hand side below is finite. Then
\begin{multline*}
\int_{\R^N}|\nabla W|^2|\xi|^2e^{-|\xi|^2}d\xi\ \ge\ D_\ell\int_{\R^N}|W|^2|\xi|^2e^{-|\xi|^2}d\xi\\
+4\inf_{d\in\C}\int_{\R^N}\big|W-d|\xi|^{\frac{D_\ell}{2}}Y(\xi/|\xi|)\big|^2|\xi|^2e^{-|\xi|^2}d\xi ,
\end{multline*}
with equality if and only if $w(\rho)=\rho^{D_\ell/2}(a_0+a_1\rho^2)$. Moreover, if $w(\rho)=\rho^{D_\ell/2}\big(C_T(N,\ell)-\rho^2\big)$, then $W$ is orthogonal to $|\xi|^{D_\ell/2}Y(\xi/|\xi|)$ in $L^2(|\xi|^2e^{-|\xi|^2}d\xi)$, and $\int|\nabla W|^2|\xi|^2e^{-|\xi|^2}=(D_\ell+4)\int|W|^2|\xi|^2e^{-|\xi|^2}$.
\end{lemma}

\begin{proof}
Let $s=D_\ell/2$, so that $s^2+Ns=c_\ell$, and write $w=\rho^sv$. Since $\int_{\bS^{N-1}}|\nabla_\sigma Y|^2d\sigma=c_\ell\int_{\bS^{N-1}}|Y|^2d\sigma$, we have
\[
\int_{\R^N}|\nabla W|^2|\xi|^2e^{-|\xi|^2}d\xi=\int_{\bS^{N-1}}|Y|^2d\sigma\int_0^\infty\Big(|w'|^2+\frac{c_\ell}{\rho^2}|w|^2\Big)\rho^{N+1}e^{-\rho^2}d\rho .
\]
Now $|w'|^2=\rho^{2s}\big(|v'|^2+\frac s\rho(|v|^2)'+\frac{s^2}{\rho^2}|v|^2\big)$, and
\[
\int_0^\infty\frac s\rho(|v|^2)'\rho^{2s+N+1}e^{-\rho^2}d\rho=-s\int_0^\infty|v|^2\big((2s+N)\rho^{2s+N-1}-2\rho^{2s+N+1}\big)e^{-\rho^2}d\rho .
\]
The coefficient of $\int|v|^2\rho^{2s+N-1}e^{-\rho^2}$ is $s^2+c_\ell-s(2s+N)=0$. Since $2s+N+2=2C_T(N,\ell)$, we get
\begin{multline}\label{subst}
\int_0^\infty\Big(|w'|^2+\frac{c_\ell}{\rho^2}|w|^2\Big)\rho^{N+1}e^{-\rho^2}d\rho\\
=\int_0^\infty|v'|^2\rho^{2C_T(N,\ell)-1}e^{-\rho^2}d\rho+D_\ell\int_0^\infty|v|^2\rho^{2C_T(N,\ell)-1}e^{-\rho^2}d\rho ,
\end{multline}
and $\int_0^\infty|w|^2\rho^{N+1}e^{-\rho^2}d\rho=\int_0^\infty|v|^2\rho^{2C_T(N,\ell)-1}e^{-\rho^2}d\rho$. The lemma follows from Theorem \ref{TA} with $\alpha=2C_T(N,\ell)$, whose second level is spanned by $C_T(N,\ell)-\rho^2$.
\end{proof}

\subsection{The Fourier transform}
We use the Fourier transform $\cF h(\xi)=\int_{\R^m}e^{-iX\cdot\xi}h(X)dX$ on $\R^m$, for which Plancherel's identity reads $\int_{\R^m}|h|^2dX=(2\pi)^{-m}\int_{\R^m}|\cF h|^2d\xi$. When the dimension has to be specified we write $\cF_m$. The Fourier transform of a radial function is radial, and for $a>-m$ one has
\begin{equation}\label{FT}
\cF_m^{-1}\big(|\cdot|^ae^{-\frac{|\cdot|^2}{2}}\big)(X)=\frac{2^{\frac a2}}{(2\pi)^{\frac m2}}\frac{\Gamma(\frac{m+a}{2})}{\Gamma(\frac m2)}\,{}_1F_1\Big(\frac{m+a}{2},\ \frac m2,\ -\frac{|X|^2}{2}\Big),
\end{equation}
see \cite[Section 3.4]{DDLLNsol}, where this formula is derived from the Hankel transform and Kummer's transformation ${}_1F_1(p,q,z)=e^z{}_1F_1(q-p,q,-z)$. We also need the Bochner--Hecke formula \cite[Chapter IV, Theorem 3.10]{SW71}. If $P$ is a homogeneous harmonic polynomial of degree $j$ on $\R^m$ and $\varphi$ is a radial profile, then
\begin{equation}\label{BH}
\cF_m\big(P\,\varphi(|\cdot|)\big)(\xi)=(-i)^j(2\pi)^{-j}\,P(\xi)\,\cF_{m+2j}\big(\varphi(|\cdot|)\big)(\xi),
\end{equation}
where the factor $(2\pi)^{-j}$ comes from our normalization of the Fourier transform. Combining \eqref{FT} and \eqref{BH} we obtain the following explicit formula, which produces all the extremals of this paper.

\begin{lemma}\label{lem:FT}
Let $m\ge1$, $j\ge0$, $a>-m$, $\lambda>0$, and let $\Psi$ be a spherical harmonic of degree $j$ on $\bS^{m-1}$. Then
\begin{multline*}
\cF_m^{-1}\Big(|\xi|^ae^{-\lambda\frac{|\xi|^2}{2}}\Psi\Big(\frac{\xi}{|\xi|}\Big)\Big)(x)\\
=\frac{i^j\,2^{\frac{a-m-j}{2}}}{\pi^{\frac m2}\lambda^{\frac{m+a+j}{2}}}\frac{\Gamma\big(\frac{m+a+j}{2}\big)}{\Gamma\big(\frac{m+2j}{2}\big)}\,{}_1F_1\Big(\frac{m+a+j}{2},\ \frac{m+2j}{2},\ -\frac{|x|^2}{2\lambda}\Big)|x|^j\Psi\Big(\frac{x}{|x|}\Big).
\end{multline*}
\end{lemma}

\begin{proof}
By a dilation it is enough to consider $\lambda=1$. Let $P(\xi)=|\xi|^j\Psi(\xi/|\xi|)$, so that the function to be transformed is $P(\xi)\,|\xi|^{a-j}e^{-|\xi|^2/2}$. By \eqref{BH}, $\cF_m^{-1}(P\varphi(|\cdot|))=i^j(2\pi)^jP\,\cF_{m+2j}^{-1}(\varphi(|\cdot|))$ for every radial profile $\varphi$. Then \eqref{FT} in dimension $m+2j$, with the exponent $a-j>-(m+2j)$, gives the formula, since $(2\pi)^j\,2^{\frac{a-j}{2}}(2\pi)^{-\frac{m+2j}{2}}=2^{\frac{a-m-j}{2}}\pi^{-\frac m2}$.
\end{proof}

\subsection{The poloidal--toroidal decomposition}\label{sec:sol}

We recall the setting of \cite{Ham23}, see also \cite{Wec97} and \cite[Section 3.1]{DDLLNsol}. Let $\partial:=r\partial_r$ with $\partial_r=\sigma\cdot\nabla$, let $\partial':=\partial+N-1$, and let $\nabla_\sigma$ and $\Delta_\sigma$ be the gradient and the Laplace--Beltrami operator on $\bS^{N-1}$. A vector field $\bU$ is \emph{toroidal} if $\bU\cdot x=0$ and $\dv\bU=0$, and \emph{poloidal} if it is of the form
\[
\bU=\bD f:=\sigma\Delta_\sigma f-\nabla_\sigma(\partial'f)
\]
for a scalar potential $f$ with vanishing spherical means. Such a field is solenoidal, see Lemma \ref{lem:lap}(ii) below. The basic facts are the following.
\begin{itemize}
\item[(a)] Every solenoidal field decomposes uniquely as $\bU=\bU_P+\bU_T$, with $\bU_P$ poloidal and $\bU_T$ toroidal.
\item[(b)] Expanding the potential in spherical harmonics, $f=\sum_{k\ge1}\sum_jf_{k,j}(r)\phi_{k,j}(\sigma)$, where $\{\phi_{k,j}\}_j$ is an orthogonal basis of $\cH_k$ such that $\int_{\bS^{N-1}}|\phi_{k,j}|^2d\sigma=|\bS^{N-1}|$, one has $\bU_P=\sum_{k\ge1}\bU_{P,k}$ with $\bU_{P,k}=\sum_j\bD(f_{k,j}\phi_{k,j})$.
\item[(c)] On every sphere, the fields of the form $a(r)\phi(\sigma)\sigma+b(r)\nabla_\sigma\phi(\sigma)$ with $\phi=\phi_{k,j}$ are orthogonal in $L^2(\bS^{N-1})$ to the fields of the same form with $\phi=\phi_{k',j'}\ne\phi_{k,j}$ and to the toroidal fields. The same holds for their gradients.
\end{itemize}
For $N\ge3$ these facts are proved in \cite[Section 2]{Ham23}. For $N=2$ they are easy to check. A solenoidal field is then $\bU=\nabla^\perp\psi$, with $\nabla^\perp=(-\partial_2,\partial_1)$. Its toroidal part is $\nabla^\perp$ of the spherical mean of $\psi$, and the toroidal fields are exactly the fields $h(|x|)x^\perp$: indeed $\bU\cdot x=0$ forces $\bU=q\,x^\perp$, and then $\dv\bU=x^\perp\cdot\nabla q=0$ forces $q$ to be radial. The facts (b) and (c) then follow from the orthogonality of the functions $e^{ik\theta}$ on the circle. We call $\bU_{P,k}$ the poloidal part of degree $k$. Since $\Delta_\sigma\phi_{k,j}=-c_k\phi_{k,j}$, one has
\begin{equation}\label{Dform}
\bD(f_k\phi_k)=-c_kf_k\phi_k\sigma-(\partial'f_k)\nabla_\sigma\phi_k\qquad(\phi_k\in\cH_k),
\end{equation}
so that the poloidal fields have the form in (c). Finally, a toroidal field has vanishing spherical means, $\int_{\bS^{N-1}}\bU(r\sigma)d\sigma=0$ for all $r>0$, see \cite[Section 3.1]{DDLLNsol} and Subsection \ref{sec:torsph}.

\begin{lemma}\label{lem:lap}\leavevmode
\begin{itemize}
\item[(i)] If $\bU$ is toroidal, then $\Delta\bU$ is toroidal.
\item[(ii)] For every scalar function $f$, $\bD f=\bm{L}(rf)$, where $\bm{L}h:=x\Delta h-\nabla\big(x\cdot\nabla h+(N-2)h\big)$. In particular, $\dv\bD f=0$.
\end{itemize}
\end{lemma}

\begin{proof}
(i) One has $\Delta(x\cdot\bU)=x\cdot\Delta\bU+2\dv\bU$ and $\dv\Delta\bU=\Delta\dv\bU$, so that $x\cdot\Delta\bU=0$ and $\dv\Delta\bU=0$.

(ii) We have $\Delta=r^{-2}\big(\partial^2+(N-2)\partial+\Delta_\sigma\big)$ and $\partial(rh)=r(\partial+1)h$. Then $x\cdot\nabla(rf)+(N-2)rf=r(\partial+N-1)f$. Hence the spherical part of $\bm{L}(rf)$ is $-r^{-1}\nabla_\sigma\big(r(\partial+N-1)f\big)=-\nabla_\sigma(\partial'f)$, and its radial component is
\[
r\Delta(rf)-\partial_r\big(r(\partial+N-1)f\big)=\Big[(\partial+1)^2+(N-2)(\partial+1)+\Delta_\sigma\Big]f-(\partial+1)(\partial+N-1)f=\Delta_\sigma f .
\]
This is $\bD f$. Finally, $\dv(x\Delta h)=N\Delta h+x\cdot\nabla\Delta h=\Delta\big(x\cdot\nabla h+(N-2)h\big)$, so that $\dv\bm{L}h=0$ for every $h$.
\end{proof}

By Lemma \ref{lem:lap}(ii), $\bD(r^{k-1}p_k\phi)=\bm{L}\big(p_k(|x|)\phi(x)\big)$ for $\phi\in\cH_k$, so that the fields of $\cE_k$ are smooth. Let $\bU=\bU_T+\sum_k\bU_{P,k}$ be a solenoidal field. By Lemma \ref{lem:lap}(i), $\Delta\bU_T$ is toroidal, and by \eqref{rLap} below, $\Delta\bD(f_{k,j}\phi_{k,j})$ has the form in (c) with $\phi=\phi_{k,j}$. Hence, by (c), each of the three integrals in \eqref{sumdef} splits, and
\begin{gather}
\delta_{2,N}(\bU)=\delta_{2,N}(\bU_T)+\sum_{k\ge1}\sum_j\delta_{2,N}\big(\bD(f_{k,j}\phi_{k,j})\big),\label{split}\\
\|\bU\|^2=\|\bU_T\|^2+\sum_{k\ge1}\sum_j\big\|\bD(f_{k,j}\phi_{k,j})\big\|^2 .\notag
\end{gather}
The same orthogonality holds for the inner product of $\|\cdot\|$, that is, these fields are mutually orthogonal in $\dot H^1(\R^N)$.

\section{The poloidal fields}\label{sec:pol}

\subsection{The three integrals}
In this section, except in the last paragraph, we fix $k\ge1$ and $\phi_k\in\cH_k$ with $\int_{\bS^{N-1}}|\phi_k|^2d\sigma=|\bS^{N-1}|$, and we consider the poloidal field $\bU_{P,k}=\bD(f_k\phi_k)$. We write $f_k=r^{k-1}g_k$ and
\[
G_k(X):=g_k(|X|),\qquad X\in\R^{N+2k}.
\]
For a smooth solenoidal field with compact support, $g_k$ is a smooth even function with compact support, since $-c_kr^kg_k\phi_k$ is the component of degree $k$ of the smooth function $x\cdot\bU$, and the component of degree $k$ of a smooth function is $r^k$ times a smooth even function of $r$. Hence $G_k$ is smooth with compact support on $\R^{N+2k}$, and the computations below are justified.

\begin{proposition}\label{prop:lift}
We have
\begin{gather*}
\int_{\R^N}|\nabla\bU_{P,k}|^2dx=\frac{c_k|\bS^{N-1}|}{|\bS^{N+2k-1}|}\int_{\R^{N+2k}}|\Delta G_k|^2dX,\\
\int_{\R^N}|\Delta\bU_{P,k}|^2dx=\frac{c_k|\bS^{N-1}|}{|\bS^{N+2k-1}|}\int_{\R^{N+2k}}|\nabla\Delta G_k|^2dX,\\
\int_{\R^N}|\nabla\bU_{P,k}|^2|x|^2dx=\frac{c_k|\bS^{N-1}|}{|\bS^{N+2k-1}|}\Big[\int_{\R^{N+2k}}|\Delta G_k|^2|X|^2dX+2\int_{\R^{N+2k}}|\nabla G_k|^2dX\Big].
\end{gather*}
\end{proposition}

\begin{proof}
We write $\bU=\bU_{P,k}$, $f=f_k$, $g=g_k$ and $\phi=\phi_k$. We use $\sigma\cdot\nabla_\sigma\phi=0$, $\int_{\bS^{N-1}}|\nabla_\sigma\phi|^2d\sigma=c_k|\bS^{N-1}|$, and the fact that $\partial$ does not act on $\sigma$ and $\nabla_\sigma\phi$. By \eqref{Dform},
\[
\partial\bU=-c_k(\partial f)\phi\sigma-(\partial'\partial f)\nabla_\sigma\phi ,
\]
and, with $\Delta_\sigma$ acting on each component,
\[
-\Delta_\sigma\bU=c_k\big(2\partial f-(c_k-N+1)f\big)\phi\sigma+\big(2c_kf-(c_k-N+3)\partial'f\big)\nabla_\sigma\phi .
\]
Since $|\nabla\bU|^2=r^{-2}\big(|\partial\bU|^2+\sum_i|\nabla_\sigma U_i|^2\big)$ and $\int_{\bS^{N-1}}|\nabla_\sigma U_i|^2d\sigma=\int_{\bS^{N-1}}(-\Delta_\sigma U_i)U_id\sigma$, we get, for $\gamma\in\{0,1\}$,
\begin{align*}
\frac{1}{c_k|\bS^{N-1}|}\int_{\R^N}|\nabla\bU|^2|x|^{2\gamma}dx=\int_0^\infty\Big[&c_k|\partial f|^2+|\partial'\partial f|^2-c_k\big(2\partial f-(c_k-N+1)f\big)f\\
&-\partial'f\big(2c_kf-(c_k-N+3)\partial'f\big)\Big]r^{N+2\gamma-3}dr .
\end{align*}
We expand the squares and integrate by parts, using
\[
2\int_0^\infty h\,\partial h\,r^adr=-(a+1)\int_0^\infty |h|^2r^adr .
\]
 The boundary terms vanish, since $f=r^{k-1}g$ with $g$ smooth, even and with compact support. In the only case where the integral of $|f|^2r^{N+2\gamma-3}$ could diverge at $0$, namely $N=2$, $k=1$ and $\gamma=0$, the terms in $|f|^2$ and in $f\,\partial f$ cancel before the integration by parts, since $c_1-N+1=0$. This gives
\begin{align*}
&\frac{1}{c_k|\bS^{N-1}|}\int_{\R^N}|\nabla\bU|^2|x|^{2\gamma}dx\\
&=\int_0^\infty\Big[|\partial^2f|^2+\big(2c_k+2-2\gamma(N-1)\big)|\partial f|^2\\
&\qquad\qquad+(c_k-N+1)\big(c_k+N-3-2\gamma(N-3)\big)|f|^2\Big]r^{N+2\gamma-3}dr .
\end{align*}
We now put $f=r^{k-1}g$, so that $\partial f=r^{k-1}\big(\partial g+(k-1)g\big)$ and $\partial^2f=r^{k-1}\big(\partial^2g+2(k-1)\partial g+(k-1)^2g\big)$. Integrating by parts again, we obtain
\begin{align}
\frac{1}{c_k|\bS^{N-1}|}\int_{\R^N}|\nabla\bU|^2dx&=\int_0^\infty\Big(|\partial^2g|^2+2(N+2k-2)|\partial g|^2\Big)r^{N+2k-5}dr,\label{g1}\\
\frac{1}{c_k|\bS^{N-1}|}\int_{\R^N}|\nabla\bU|^2|x|^2dx&=\int_0^\infty\Big(|\partial^2g|^2+2|\partial g|^2\Big)r^{N+2k-3}dr.\label{g2}
\end{align}
Next we compute $\Delta\bU$. Since $r^2\Delta=\partial^2+(N-2)\partial+\Delta_\sigma$ on each component, \eqref{Dform} and the formula for $\Delta_\sigma\bU$ give, after the substitution $f=r^{k-1}g$,
\begin{equation}\label{rLap}
r^2\Delta\bU=-c_kr^{k-1}\big(\partial^2g+(N+2k-2)\partial g\big)\phi\sigma-r^{k-1}\big(\partial^3g+A_k\partial^2g+B_k\partial g\big)\nabla_\sigma\phi ,
\end{equation}
where $A_k=2N+3k-6$ and $B_k=2k^2+3kN-10k+(N-2)(N-4)$. Hence
\begin{multline*}
\frac{1}{c_k|\bS^{N-1}|}\int_{\R^N}|\Delta\bU|^2dx\\
=\int_0^\infty\Big[c_k\big|\partial^2g+(N+2k-2)\partial g\big|^2+\big|\partial^3g+A_k\partial^2g+B_k\partial g\big|^2\Big]r^{N+2k-7}dr .
\end{multline*}
We expand the squares and integrate by parts as above. Since, with $m=N+2k$,
\begin{gather*}
c_k+A_k^2-A_k(m-6)-2B_k=6m-16,\\ 4c_k(m-2)+B_k^2-A_kB_k(m-6)+B_k(m-6)^2=8(m-2)(m-4),
\end{gather*}
we arrive at
\begin{multline}\label{g3}
\frac{1}{c_k|\bS^{N-1}|}\int_{\R^N}|\Delta\bU|^2dx\\
=\int_0^\infty\Big(|\partial^3g|^2+\big(6(N+2k)-16\big)|\partial^2g|^2+8(N+2k-2)(N+2k-4)|\partial g|^2\Big)r^{N+2k-7}dr .
\end{multline}
Finally, for a radial function $H(X)=h(|X|)$ on $\R^m$, the same computation gives
\begin{gather*}
\int_{\R^m}|\Delta H|^2dX=|\bS^{m-1}|\int_0^\infty\Big(|\partial^2h|^2+2(m-2)|\partial h|^2\Big)r^{m-5}dr,\\
\int_{\R^m}|\Delta H|^2|X|^2dX=|\bS^{m-1}|\int_0^\infty|\partial^2h|^2r^{m-3}dr,\\
\int_{\R^m}|\nabla H|^2dX=|\bS^{m-1}|\int_0^\infty|\partial h|^2r^{m-3}dr,
\end{gather*}
and
\begin{multline*}
\int_{\R^m}|\nabla\Delta H|^2dX\\
=|\bS^{m-1}|\int_0^\infty\Big(|\partial^3h|^2+(6m-16)|\partial^2h|^2+8(m-2)(m-4)|\partial h|^2\Big)r^{m-7}dr .
\end{multline*}
We apply these formulas with $m=N+2k$ and $H=G_k$, and we compare them with \eqref{g1}, \eqref{g2} and \eqref{g3}. This proves the proposition.
\end{proof}

Since both sides of the identities in Proposition \ref{prop:lift} are quadratic forms in $g_k$, the same identities hold for the associated bilinear forms.

\subsection{The Fourier transform}
The Fourier transform of $G_k$ on $\R^{N+2k}$ is radial and real-valued. We write it as
\[
\cF G_k(\xi)=|\xi|^{\beta_k}F_k(\xi),\qquad\xi\in\R^{N+2k},
\]
with $\beta_k$ given by \eqref{betaD}. In this subsection and in the next one, all the integrals in $\xi$ are over $\R^{N+2k}$. By Proposition \ref{prop:lift}, Plancherel's identity and $\cF(X_j\Delta G_k)=-i\partial_{\xi_j}(|\xi|^2\cF G_k)$,
\begin{multline*}
\int_{\R^N}|\Delta\bU_{P,k}|^2dx+\int_{\R^N}|\nabla\bU_{P,k}|^2|x|^2dx=\frac{c_k|\bS^{N-1}|}{(2\pi)^{N+2k}|\bS^{N+2k-1}|}\Big[\int|\cF G_k|^2|\xi|^6d\xi\\
+\int\big|\nabla\big(|\xi|^2\cF G_k\big)\big|^2d\xi+2\int|\cF G_k|^2|\xi|^2d\xi\Big].
\end{multline*}
An integration by parts of the cross term gives
\[
\int\big|\nabla\big(|\xi|^2\cF G_k\big)\big|^2d\xi=\int|\nabla F_k|^2|\xi|^{4+2\beta_k}d\xi-(2+\beta_k)(N+2k+\beta_k)\int|F_k|^2|\xi|^{2+2\beta_k}d\xi .
\]
Since $\beta_k$ is the root of $(2+\beta)(N+2k+\beta)=2$ with $N+2k+\beta>0$, the last term cancels the term $2\int|\cF G_k|^2|\xi|^2$. Therefore
\begin{multline}\label{blockw}
\delta_{2,N}(\bU_{P,k})-2\big(C_P(N,k)-C_N\big)\|\bU_{P,k}\|^2\\
=\frac{c_k|\bS^{N-1}|}{(2\pi)^{N+2k}|\bS^{N+2k-1}|}\Big[\int|F_k|^2|\xi|^{6+2\beta_k}d\xi+\int|\nabla F_k|^2|\xi|^{4+2\beta_k}d\xi\\
-2C_P(N,k)\int|F_k|^2|\xi|^{4+2\beta_k}d\xi\Big],
\end{multline}
and, with the same constant in front,
\begin{equation}\label{normw}
\|\bU_{P,k}\|^2=\frac{c_k|\bS^{N-1}|}{(2\pi)^{N+2k}|\bS^{N+2k-1}|}\int|F_k|^2|\xi|^{4+2\beta_k}d\xi .
\end{equation}
So, after the Fourier transform, the left-hand side of \eqref{blockw} is the deficit of a first order Caffarelli--Kohn--Nirenberg inequality for the radial function $F_k$ on $\R^{N+2k}$.

\subsection{The block estimate}
By Lemma \ref{lem:FT} with $m=N+2k$, $j=0$, $a=\beta_k$ and $\lambda=1$, the function $\cF^{-1}\big(e^{-|\xi|^2/2}|\xi|^{\beta_k}\big)$ is a positive multiple of $p_k(|X|)$. Thus $F_k$ is a multiple of $e^{-|\xi|^2/2}$ if and only if $g_k$ is a multiple of $p_k$, i.e.\ $\bU_{P,k}$ is a multiple of $\bD(r^{k-1}p_k\phi_k)$. These fields have finite integrals in \eqref{2ndHUP}, since $F_k=e^{-|\xi|^2/2}$ gives finite integrals in \eqref{blockw} and \eqref{normw}.

\begin{proposition}[Poloidal block estimate]\label{prop:pol}
We have
\begin{equation}\label{polblock}
\delta_{2,N}(\bU_{P,k})\ \ge\ 4\inf_{c\in\R}\big\|\bU_{P,k}-c\,\bD(r^{k-1}p_k\phi_k)\big\|^2+2\big(C_P(N,k)-C_N\big)\|\bU_{P,k}\|^2 ,
\end{equation}
and the constant $4$ is sharp. Moreover, the left-hand side of \eqref{blockw} is nonnegative, and it vanishes if and only if $\bU_{P,k}$ is a multiple of $\bD(r^{k-1}p_k\phi_k)$.
\end{proposition}

\begin{proof}
Since $N+2k+(4+2\beta_k)=2C_P(N,k)$, we may apply Corollary \ref{cor:CKN} with $m=N+2k$ and $p=4+2\beta_k$ to the right-hand side of \eqref{blockw}. It is bounded from below by $4$ times
\[
\frac{c_k|\bS^{N-1}|}{(2\pi)^{N+2k}|\bS^{N+2k-1}|}\inf_{c\in\C}\int\big|F_k-ce^{-\frac{|\xi|^2}{2}}\big|^2|\xi|^{4+2\beta_k}d\xi .
\]
The infimum may be taken over real numbers, since $F_k$ is real-valued. Since $|F_k-ce^{-|\xi|^2/2}|^2|\xi|^{4+2\beta_k}=\big||\xi|^2\big(\cF G_k-c|\xi|^{\beta_k}e^{-|\xi|^2/2}\big)\big|^2$, Plancherel's identity and the first identity of Proposition \ref{prop:lift} show that this quantity is $\inf_{c\in\R}\|\bU_{P,k}-c\,\bD(r^{k-1}p_k\phi_k)\|^2$. This proves \eqref{polblock}. The other statements follow from the cases of equality in Corollary \ref{cor:CKN}. In particular, the constant $4$ is attained when $F_k=(|\xi|^2-C_P(N,k))e^{-|\xi|^2/2}$.
\end{proof}

By Corollary \ref{cor:CKN}, equality holds in \eqref{polblock} if and only if $F_k=(a_0+a_1|\xi|^2)e^{-|\xi|^2/2}$. Let $\phi\in\cH_k\setminus\{0\}$, and let $\bR$ be the poloidal field $\bD(f\phi)$ whose function $F_k$ is $\big(C_P(N,k)-|\xi|^2\big)e^{-|\xi|^2/2}$. By Lemma \ref{lem:FT}, its profile is a linear combination of two Kummer functions, and it is smooth with finite integrals. By Theorem \ref{TA} with $\alpha=2C_P(N,k)$, the functions $e^{-|\xi|^2/2}$ and $\big(C_P(N,k)-|\xi|^2\big)e^{-|\xi|^2/2}$ are orthogonal in $L^2(|\xi|^{4+2\beta_k}d\xi)$. By the polarization of \eqref{normw}, $\bR$ is therefore orthogonal in $\dot H^1$ to $\bD(r^{k-1}p_k\phi)$, and \eqref{blockw} gives
\begin{equation}\label{Rpol}
\delta_{2,N}(\bR)=\big(4+2C_P(N,k)-2C_N\big)\|\bR\|^2 .
\end{equation}
We call $\bR$ a radial excitation of $\bD(r^{k-1}p_k\phi)$. Indeed, the derivative at $\lambda=1$ of the dilations of $\bD(r^{k-1}p_k\phi)$ is a linear combination of $\bD(r^{k-1}p_k\phi)$ and $\bR$. To see this, note that the dilations act on $\cF G_k$ by $\lambda\mapsto|\xi/\lambda|^{\beta_k}e^{-|\xi|^2/(2\lambda^2)}$, up to a power of $\lambda$ coming from the normalization in \eqref{dilation}, and that the derivative of this function at $\lambda=1$ is $(|\xi|^2-\beta_k)|\xi|^{\beta_k}e^{-|\xi|^2/2}$. The power of $\lambda$ only adds a multiple of $\bD(r^{k-1}p_k\phi)$.

Let now $\bU_{P,k}=\sum_j\bD(f_{k,j}\phi_{k,j})$ be the poloidal part of degree $k$ of a solenoidal field. We apply \eqref{polblock} to each field $\bD(f_{k,j}\phi_{k,j})$ and we sum over $j$. Since the fields $\bD(r^{k-1}p_k\phi_{k,j})$ span $\{\bD(r^{k-1}p_k\phi):\phi\in\cH_k\}$ and are mutually orthogonal in $\dot H^1$, we obtain
\begin{equation}\label{polk}
\delta_{2,N}(\bU_{P,k})\ \ge\ 4\dist^2\big(\bU_{P,k},\{\bD(r^{k-1}p_k\phi):\phi\in\cH_k\}\big)+2\big(C_P(N,k)-C_N\big)\|\bU_{P,k}\|^2 .
\end{equation}

\section{The toroidal fields}\label{sec:tor}

\subsection{Toroidal fields and spherical harmonics}\label{sec:torsph}
We use polar coordinates. For a vector field $\bF$ we write $\bF=F_r\sigma+\bF_{\rm tan}$ with $F_r=\sigma\cdot\bF$. Then
\begin{equation}\label{divpolar}
\dv\bF=\partial_rF_r+\frac{N-1}{r}F_r+\frac1r\dv_\sigma\bF_{\rm tan},
\end{equation}
where $\dv_\sigma$ is the divergence on $\bS^{N-1}$ applied to $\sigma\mapsto\bF_{\rm tan}(r\sigma)$. A field $\bZ$ on $\bS^{N-1}$ is tangent if $\sigma\cdot\bZ(\sigma)=0$ for all $\sigma$. If $\bU$ is toroidal, then for every $r>0$ the field $\sigma\mapsto\bU(r\sigma)$ is tangent, and it is divergence-free on $\bS^{N-1}$ by \eqref{divpolar}. Conversely, let $\bP\in\cH_\ell^N$ be such that its restriction to $\bS^{N-1}$ is tangent and divergence-free on $\bS^{N-1}$. Then $x\cdot\bP(x)=r^{\ell+1}\sigma\cdot\bP(\sigma)=0$, and \eqref{divpolar} gives $\dv\bP(r\sigma)=r^{\ell-1}\dv_\sigma\bP(\sigma)=0$. So $\bP$ is toroidal. Also, $h(|x|)\bP$ is toroidal for every radial profile $h$, since $x\cdot(h\bP)=0$ and $\dv(h\bP)=h'\frac{x\cdot\bP}{|x|}+h\dv\bP=0$.

For a vector field $\bV$ on $\R^N$ and $\ell\ge0$, let $\bV^{(\ell)}$ be the field obtained by projecting each component of $\bV$ onto $\cH_\ell$ on every sphere. We say that $\bV$ has degree $\ell$ if $\bV=\bV^{(\ell)}$. The three integrals in \eqref{sumdef} are sums over the components of scalar quadratic forms which commute with rotations. Hence they split over the degrees,
\begin{equation}\label{splitT}
\delta_{2,N}(\bV)=\sum_{\ell\ge0}\delta_{2,N}(\bV^{(\ell)}),\qquad \|\bV\|^2=\sum_{\ell\ge0}\|\bV^{(\ell)}\|^2 ,
\end{equation}
and the parts of different degrees are orthogonal in $\dot H^1(\R^N)$. This holds for all vector fields, toroidal or not. The same notation is used on $\bS^{N-1}$.

Let $\bZ$ be a smooth tangent field on $\bS^{N-1}$ with $\dv_\sigma\bZ=0$. Then $\bZ^{(0)}=0$. Indeed, let $\mathbf e_i$ be the $i$-th unit vector. The tangential gradient of $\sigma_i$ is $\nabla_\sigma\sigma_i=\mathbf e_i-\sigma_i\sigma$, so that $\bZ\cdot\nabla_\sigma\sigma_i=\bZ\cdot\mathbf e_i$, and
\[
\int_{\bS^{N-1}}\bZ\cdot\mathbf e_i\,d\sigma=\int_{\bS^{N-1}}\bZ\cdot\nabla_\sigma\sigma_i\,d\sigma=-\int_{\bS^{N-1}}\sigma_i\dv_\sigma\bZ\,d\sigma=0 .
\]
Applying this to $\sigma\mapsto\bU_T(r\sigma)$ for every $r>0$, we get $\bU_T^{(0)}=0$ for every toroidal field $\bU_T$. For $N=2$ moreover $\bU_T=h(|x|)x^\perp=\bU_T^{(1)}$.

\subsection{The Fourier transform}
Let $\bU_T$ be a toroidal field, and let $\bW_T:=e^{|\xi|^2/2}\cF\bU_T$. Since $\dv\bU_T=0$ and $x\cdot\bU_T=0$, we have $\xi\cdot\cF\bU_T=0$ and $\dv\cF\bU_T=0$, so that $\cF\bU_T$ and $\bW_T$ are toroidal. By Plancherel's identity,
\[
\int|\Delta\bU_T|^2dx=\frac{1}{(2\pi)^N}\int|\xi|^4|\cF\bU_T|^2d\xi,\qquad
\|\bU_T\|^2=\frac{1}{(2\pi)^N}\int|\xi|^2|\cF\bU_T|^2d\xi,\]
\[\int|\nabla\bU_T|^2|x|^2dx=\frac{1}{(2\pi)^N}\int|\xi|^2|\nabla\cF\bU_T|^2d\xi .
\]
For the last identity, $\cF(x_j\partial_iU_n)=-\partial_{\xi_j}(\xi_i\cF U_n)$, where $U_n$ are the components of $\bU_T$, and the sum of the first two terms on the right-hand side below vanishes:
\[
\sum_{i,j,n}\int\big|\delta_{ij}\cF U_n+\xi_i\partial_j\cF U_n\big|^2d\xi=N\int|\cF\bU_T|^2d\xi+\int\xi\cdot\nabla\big(|\cF\bU_T|^2\big)d\xi+\int|\xi|^2|\nabla\cF\bU_T|^2d\xi .
\]
By \eqref{gauss} with $m=N$ and $p=2$, we obtain
\begin{equation}\label{torw}
\delta_{2,N}(\bU_T)=\frac{1}{(2\pi)^N}\Big[\int_{\R^N}|\nabla\bW_T|^2|\xi|^2e^{-|\xi|^2}d\xi-(2C_N-N-2)\int_{\R^N}|\bW_T|^2|\xi|^2e^{-|\xi|^2}d\xi\Big],
\end{equation}
and
\begin{equation}\label{normv}
\|\bU_T\|^2=\frac{1}{(2\pi)^N}\int_{\R^N}|\bW_T|^2|\xi|^2e^{-|\xi|^2}d\xi .
\end{equation}
The proof of \eqref{torw} and \eqref{normv} uses only Plancherel's identity and \eqref{gauss}, so these identities hold for every vector field. Since the Fourier transform commutes with rotations, the part of degree $\ell$ of $\bW_T$ is $e^{|\xi|^2/2}\cF\big(\bU_T^{(\ell)}\big)$. We denote it by $\bW_T^{(\ell)}$. Hence \eqref{torw} and \eqref{normv} hold for each $\bU_T^{(\ell)}$ and $\bW_T^{(\ell)}$, and $\bW_T^{(0)}=0$.

\subsection{The block estimate}
Let $\ell\ge1$, and let $\{Y_j\}_j$ be an orthonormal basis of $\cH_\ell$ in $L^2(\bS^{N-1})$. We write the components of $\bW_T^{(\ell)}$ as $\sum_jw_{ij}(\rho)Y_j(\sigma)$. We apply Lemma \ref{lem:degree} to each function $w_{ij}(\rho)Y_j(\sigma)$, and we sum. Then \eqref{torw}, \eqref{normv} and the identity $2C_T(N,\ell)=N+2+D_\ell$ give
\begin{multline}\label{torblockW}
\delta_{2,N}(\bU_T^{(\ell)})\ \ge\ 2\big(C_T(N,\ell)-C_N\big)\|\bU_T^{(\ell)}\|^2\\
+\frac{4}{(2\pi)^N}\inf_{\bP}\int_{\R^N}\big|\bW_T^{(\ell)}-|\xi|^{\frac{D_\ell}{2}-\ell}\bP(\xi)\big|^2|\xi|^2e^{-|\xi|^2}d\xi ,
\end{multline}
where the infimum is taken over all $\bP\in\cH_\ell^N$ with complex coefficients. Equality holds if and only if all the functions $w_{ij}$ are of the form $\rho^{D_\ell/2}(a_0+a_1\rho^2)$. Here we do not use that $\bU_T^{(\ell)}$ is toroidal. The toroidal condition on $\bU_T$ enters only through the vanishing of $\bW_T^{(0)}$ and through Lemma \ref{lem:proj} below, which shows that the best approximation in \eqref{torblockW} is toroidal. For $\ell=1$ this is a direct computation. For $\ell\ge2$ its proof uses the following lemma on the sphere, which we prove for the reader's convenience.

\begin{lemma}\label{lem:sphere}
Let $\bZ$ be a smooth tangent field on $\bS^{N-1}$ with $\dv_\sigma\bZ=0$. Then for every $\ell\ge1$ the field $\bZ^{(\ell)}$ is tangent and divergence-free on $\bS^{N-1}$.
\end{lemma}

\begin{proof}
Let $\bF=\sum_{\ell\ge0}\bF_\ell$, with $\bF_\ell(x)=|x|^\ell\bZ^{(\ell)}(x/|x|)$, be the harmonic extension of $\bZ$ to the unit ball, component by component. Each component of $\bF_\ell$ is a harmonic polynomial, homogeneous of degree $\ell$, and the series converges in $C^\infty$ of the closed ball since $\bZ$ is smooth. Let $g=x\cdot\bF$ and $h=\dv\bF$. Then $h$ is harmonic and $\Delta g=2h$, so that $g$ is biharmonic. On the sphere, $g=\sigma\cdot\bZ=0$. Since $F_r=g/r$, \eqref{divpolar} at $r=1$ gives $h=\partial_rg$ on the sphere. We expand $g$ in spherical harmonics. A smooth biharmonic function in the ball has the components $(ar^k+br^{k+2})Y$ with $Y\in\cH_k$. The condition $g=0$ on the sphere gives $a=-b$. Since $\Delta(r^{k+2}Y)=(4k+2N)r^kY$, the condition $\Delta g=2\partial_rg$ on the sphere gives $b(4k+2N)=4b$, that is, $b=0$, except when $k=0$ and $N=2$. In this last case all the other components of $g$ vanish, so that $h=\frac12\Delta g$ is constant. Since $\int_{\{|x|<1\}}h\,dx=\int_{\bS^{N-1}}\sigma\cdot\bZ\,d\sigma=0$, we get $b=0$ again. Therefore $g=0$ and $h=0$ in the ball. Taking the homogeneous parts of $x\cdot\bF=0$ and of $\dv\bF=0$, we get $x\cdot\bF_\ell=0$ and $\dv\bF_\ell=0$ for every $\ell$. Hence $\bZ^{(\ell)}$ is tangent, and \eqref{divpolar} at $r=1$ gives $\dv_\sigma\bZ^{(\ell)}=\dv\bF_\ell=0$.
\end{proof}

\begin{lemma}\label{lem:proj}
The infimum in \eqref{torblockW} is attained at a toroidal $\bP$.
\end{lemma}

\begin{proof}
We write $\omega_s(\xi):=|\xi|^se^{-|\xi|^2}$. The infimum in \eqref{torblockW} is a minimum, since it is taken over a finite dimensional space. Since the components of $\bW_T-\bW_T^{(\ell)}$ are orthogonal to $\cH_\ell$ on every sphere, we have for every $\bP\in\cH_\ell^N$
\begin{equation}\label{fullW}
\int_{\R^N}\big|\bW_T^{(\ell)}-|\xi|^{\frac{D_\ell}{2}-\ell}\bP\big|^2\omega_2\,d\xi=\int_{\R^N}\big|\bW_T-|\xi|^{\frac{D_\ell}{2}-\ell}\bP\big|^2\omega_2\,d\xi-\int_{\R^N}\big|\bW_T-\bW_T^{(\ell)}\big|^2\omega_2\,d\xi .
\end{equation}
So we may work with the full field $\bW_T$, which is toroidal.

We first consider $\ell=1$. In this case $\bP(\xi)=A\xi$ with $A\in\C^{N\times N}$. We write $A=A_{\rm sym}+A_{\rm skew}$, with $A_{\rm sym}=\frac{A+A^T}{2}$ and $A_{\rm skew}=\frac{A-A^T}{2}$. Then
\begin{multline*}
\int\big|\bW_T-|\xi|^{\frac{D_1}{2}-1}A\xi\big|^2\omega_2=\int\big|\bW_T-|\xi|^{\frac{D_1}{2}-1}A_{\rm skew}\xi\big|^2\omega_2+\int|A_{\rm sym}\xi|^2\omega_{D_1}\\
-2{\rm Re}\int\bW_T\cdot\overline{A_{\rm sym}\xi}\,\omega_{\frac{D_1}{2}+1}+2{\rm Re}\int A_{\rm skew}\xi\cdot\overline{A_{\rm sym}\xi}\,\omega_{D_1} .
\end{multline*}
The last two integrals vanish. For the first one, let $\varphi(\xi)=\frac12\xi\cdot A_{\rm sym}\xi$, so that $\nabla\varphi=A_{\rm sym}\xi$. Since $\bW_T$ is toroidal and $\omega_{D_1/2+1}$ is radial, we have $\dv\big(\bW_T\,\omega_{D_1/2+1}\big)=\omega_{D_1/2+1}\dv\bW_T+\frac{\omega_{D_1/2+1}'(|\xi|)}{|\xi|}\,\xi\cdot\bW_T=0$, and an integration by parts gives $\int\bW_T\cdot\overline{\nabla\varphi}\,\omega_{D_1/2+1}=0$. For the second one, we write $A_{\rm skew}=(a_{ij})$ and $A_{\rm sym}=(b_{ij})$. Since $\int_{\bS^{N-1}}\sigma_j\sigma_kd\sigma=\frac{|\bS^{N-1}|}{N}\delta_{jk}$, the integral is a multiple of $\sum_{i,j}a_{ij}\overline{b_{ij}}$, which is $0$ because $a_{ji}=-a_{ij}$ and $b_{ji}=b_{ij}$. Hence
\[
\int\big|\bW_T-|\xi|^{\frac{D_1}{2}-1}A\xi\big|^2\omega_2=\int\big|\bW_T-|\xi|^{\frac{D_1}{2}-1}A_{\rm skew}\xi\big|^2\omega_2+\int|A_{\rm sym}\xi|^2\omega_{D_1},
\]
and the best matrix $A$ is antisymmetric. Then $\bP(\xi)=A\xi$ is toroidal, since $\xi\cdot A\xi=0$ and $\dv(A\xi)=\Tr A=0$.

Let now $\ell\ge2$, let $\bP_*$ be the minimizer, and let $\bQ\in\cH_\ell^N$. Since $\bQ(\rho\sigma)=\rho^\ell\bQ(\sigma)$, \eqref{fullW} shows that the minimizer is characterized by
\[
\int_{\bS^{N-1}}\bP_*\cdot\overline{\bQ}\,d\sigma\int_0^\infty\rho^{D_\ell+N+1}e^{-\rho^2}d\rho=\int_{\bS^{N-1}}\bZ\cdot\overline{\bQ}\,d\sigma\qquad\text{for all }\bQ\in\cH_\ell^N,
\]
where
\[
\bZ(\sigma):=\int_0^\infty\bW_T(\rho\sigma)\rho^{\frac{D_\ell}{2}+N+1}e^{-\rho^2}d\rho .
\]
By the convention of Section \ref{sec:intro}, we may assume that $\bU_T$ is smooth with compact support, so that $\cF\bU_T$ and $\bW_T$ are smooth. The general case follows by density, since the minimizer depends continuously on $\bW_T^{(\ell)}$ in $L^2(\omega_2d\xi)$ and the toroidal fields form a closed subspace of $\cH_\ell^N$. The real and imaginary parts of $\bZ$ are then smooth tangent fields, and they are divergence-free on $\bS^{N-1}$, since $\sigma\mapsto\bW_T(\rho\sigma)$ has these properties for every $\rho>0$. Hence the restriction of $\bP_*$ to the sphere is a positive multiple of $\bZ^{(\ell)}$, which is tangent and divergence-free by Lemma \ref{lem:sphere}. By the discussion at the beginning of this section, $\bP_*$ is toroidal.
\end{proof}

Let $\bP\in\cH_\ell^N$ be toroidal. By Lemma \ref{lem:FT} with $m=N$, $j=\ell$, $a=D_\ell/2$ and $\lambda=1$, applied to each component, the inverse Fourier transform of $|\xi|^{\frac{D_\ell}{2}-\ell}e^{-|\xi|^2/2}\bP(\xi)$ is $i^\ell$ times a positive multiple of $t_\ell(|x|)\bP(x)$. Hence, by Plancherel's identity, the infimum in \eqref{torblockW} is $(2\pi)^N$ times the infimum of $\|\bU_T^{(\ell)}-t_\ell(|x|)\bP\|^2$ over the toroidal $\bP\in\cH_\ell^N$ with complex coefficients. Since $\bU_T^{(\ell)}$ is real-valued, $\|\bU_T^{(\ell)}-t_\ell(\bP_1+i\bP_2)\|^2=\|\bU_T^{(\ell)}-t_\ell\bP_1\|^2+\|t_\ell\bP_2\|^2$ for real $\bP_1,\bP_2$. Therefore \eqref{torblockW} and Lemma \ref{lem:proj} give
\begin{equation}\label{torl}
\delta_{2,N}(\bU_T^{(\ell)})\ \ge\ 4\dist^2\big(\bU_T^{(\ell)},\{t_\ell(|x|)\bP:\bP\in\cH_\ell^N\ \text{toroidal}\}\big)+2\big(C_T(N,\ell)-C_N\big)\|\bU_T^{(\ell)}\|^2 .
\end{equation}
In the same way, for every field $\bV$ of degree $\ell$, not necessarily toroidal, Lemma \ref{lem:degree} gives
\begin{equation}\label{torall}
\int_{\R^N}|\Delta\bV|^2dx+\int_{\R^N}|\nabla\bV|^2|x|^2dx\ \ge\ 2C_T(N,\ell)\|\bV\|^2 ,
\end{equation}
with equality if and only if $\bV=t_\ell(|x|)\bP$ for some $\bP\in\cH_\ell^N$. For $\ell=1$, the next paragraph shows that the constant $4$ in \eqref{torl} is sharp.

Let $A\ne0$ be a real antisymmetric matrix, and let $\bR$ be the field whose Fourier transform is
\[
\cF\bR(\xi)=i\,|\xi|^{\frac{D_1}{2}-1}\big(C_T(N,1)-|\xi|^2\big)e^{-\frac{|\xi|^2}{2}}A\xi .
\]
By Lemma \ref{lem:FT} with $j=1$, $\bR=h(|x|)Ax$ for a smooth real profile $h$ (the factor $i$ in $\cF\bR$ and the factor $i^j=i$ of Lemma \ref{lem:FT} multiply to $-1$), so that $\bR$ is a real toroidal field of degree one. By Lemma \ref{lem:degree}, $\bR$ is orthogonal in $\dot H^1$ to all the fields $t_1(|x|)\bP$, $\bP\in\cH_1^N$, and \eqref{torw} gives
\begin{equation}\label{Rtor}
\delta_{2,N}(\bR)=\big(4+2C_T(N,1)-2C_N\big)\|\bR\|^2 .
\end{equation}
We call $\bR$ a radial excitation of $t_1(|x|)Ax$.

Finally, for $N\ge3$ and every $\ell\ge1$ there are nonzero toroidal fields in $\cH_\ell^N$. Indeed, for $Y\in\cH_\ell$ the field $\bP:=(x_2\partial_3Y-x_3\partial_2Y,\ x_3\partial_1Y-x_1\partial_3Y,\ x_1\partial_2Y-x_2\partial_1Y,\ 0,\dots,0)$ belongs to $\cH_\ell^N$, since the operators $x_i\partial_j-x_j\partial_i$ commute with $\Delta$ and preserve homogeneity. A direct computation gives $x\cdot\bP=0$ and $\dv\bP=0$, and $\bP\ne0$ for $Y=\mathrm{Re}\,(x_1+ix_2)^\ell$. Hence $C_T(N,\ell)$ is the sharp constant of \eqref{2ndHUP} on the toroidal fields of degree $\ell$. For $\ell=1$ such fields are the fields $t_1(|x|)Ax$ with $A\ne0$ antisymmetric, which exist for all $N\ge2$.

\section{The stability of the second order HUP for solenoidal fields--Proofs of Theorem \texorpdfstring{\ref{thm:sharp}}{1.1} and Theorem \texorpdfstring{\ref{thm:sum}}{1.2}}\label{sec:proofs}

\subsection{The blocks}
Let $\bU=\bU_T+\sum_{k\ge1}\bU_{P,k}$ be a solenoidal field. We call the fields $\bU_{P,k}$, $k\ge1$, and $\bU_T^{(\ell)}$, $\ell\ge1$, the blocks of $\bU$. Recall that $\bU_T^{(0)}=0$. By \eqref{split} and \eqref{splitT}, $\delta_{2,N}(\bU)$ and $\|\bU\|^2$ are the sums of the corresponding quantities over the blocks.

For a block $\bV$ we write $C:=C_P(N,k)$ and $\cE(\bV):=\{\bD(r^{k-1}p_k\phi):\phi\in\cH_k\}$ if $\bV=\bU_{P,k}$, and $C:=C_T(N,\ell)$ and $\cE(\bV):=\{t_\ell(|x|)\bP:\bP\in\cH_\ell^N\ \text{toroidal}\}$ if $\bV=\bU_T^{(\ell)}$. We call $2(C-C_N)$ the level of the block. Then \eqref{polk} and \eqref{torl} say that
\begin{equation}\label{blockest}
\delta_{2,N}(\bV)\ \ge\ 2(C-C_N)\|\bV\|^2+4\dist^2(\bV,\cE(\bV)) .
\end{equation}
Note that $\cE_k$ is the sum of the spaces $\cE(\bU_{P,j})$, $j\le k$, and $\cT_k$ is the sum of the spaces $\cE(\bU_T^{(j)})$, $j\le k$. Every field of $\cE(\bV)$ is a poloidal field of degree $k$, respectively a toroidal field of degree $\ell$. Let $\cG$ be a sum of some of the spaces $\cE(\bV)$, and let $\cG=\cG_P+\cG_T$, where $\cG_P$ is the sum of the poloidal spaces and $\cG_T$ is the sum of the toroidal spaces. Since the poloidal fields are orthogonal to the toroidal fields in $\dot H^1$, we have $\dist^2(\bU,\cG)=\dist^2(\bU_P,\cG_P)+\dist^2(\bU_T,\cG_T)$. By \eqref{split} for the poloidal part, and by \eqref{splitT} for the toroidal part, the distance splits further:
\begin{equation}\label{distsplit}
\dist^2(\bU,\cG)=\sum\dist^2(\bV,\cE(\bV))+\sum\|\bV\|^2,
\end{equation}
where the first sum is over the blocks $\bV$ with $\cE(\bV)\subset\cG$, and the second sum is over the other blocks.

\subsection{Proof of the inequality in Theorem \ref{thm:sharp}}
Let $\bU$ be a solenoidal field. By \eqref{blockest}, we have $\delta_{2,N}(\bU)\ge\sum2(C-C_N)\|\bV\|^2\ge0$, where the sum is over the blocks $\bV$ of $\bU$. For the dilates \eqref{dilation} one has
\begin{equation}\label{scaling}
\begin{gathered}
\int_{\R^N}|\Delta\bU_\lambda|^2dx=\lambda^2\int_{\R^N}|\Delta\bU|^2dx,\qquad \|\bU_\lambda\|=\|\bU\|,\\
\int_{\R^N}|\nabla\bU_\lambda|^2|x|^2dx=\lambda^{-2}\int_{\R^N}|\nabla\bU|^2|x|^2dx .
\end{gathered}
\end{equation}
Applying $\delta_{2,N}(\bU_\lambda)\ge0$ and minimizing in $\lambda$, we obtain \eqref{2ndHUPsol}. In the same way, Proposition \ref{prop:pol} and \eqref{torall}, together with the existence of toroidal fields of every degree shown at the end of Section \ref{sec:tor}, give the statements on $C_P(N,k)$ and $C_T(N,k)$. Indeed, the inequality in sum form on a class of fields which is invariant under dilations implies the inequality in product form with the same constant. A field which gives equality in the sum form also gives equality in the product form, since $2\sqrt{I_1I_2}\le I_1+I_2$ for $I_1,I_2\ge0$. The characterization of the equality cases in \eqref{2ndHUPsol} is proved in Section \ref{sec:prod}.

\subsection{The levels below \texorpdfstring{$4$}{4}}
\begin{lemma}\label{lem:levels}
The blocks with level smaller than $4$ are the following, listed in the order of their levels.
\begin{itemize}
\item[(i)] For $N\ge5$: $\bU_{P,1}$, $\bU_T^{(1)}$, $\bU_{P,2}$, $\bU_T^{(2)}$, $\bU_{P,3}$, with the levels
\begin{gather*}
0,\qquad \sqrt{N^2+4N-4}-\sqrt{N^2+8},\qquad \sqrt{N^2+4N+12}-\sqrt{N^2+8},\\
\sqrt{N^2+8N}-\sqrt{N^2+8},\qquad \sqrt{N^2+8N+24}-\sqrt{N^2+8}.
\end{gather*}
\item[(ii)] For $N=4$: the same blocks and $\bU_T^{(3)}$, with level $2\sqrt{19}-2\sqrt6$.
\item[(iii)] For $N=3$: $\bU_{P,k}$ and $\bU_T^{(k)}$, $k\le3$. The blocks of degree $k$ have the same level, which is $0$, $\sqrt{33}-\sqrt{17}$ and $\sqrt{57}-\sqrt{17}$ for $k=1,2,3$.
\item[(iv)] For $N=2$: $\bU_T^{(1)}$, $\bU_{P,1}$, $\bU_{P,2}$, $\bU_{P,3}$, with levels $0$, $2\sqrt3-2\sqrt2$, $2\sqrt6-2\sqrt2$ and $2\sqrt{11}-2\sqrt2$.
\end{itemize}
All the other blocks have levels larger than $4$.
\end{lemma}

\begin{proof}
For a poloidal block of degree $k$ let $\gamma:=c_k-N+3$, and for a toroidal block of degree $\ell$ let $\gamma:=c_\ell$. Since $(N+2k-2)^2+8=N^2+4(c_k-N+3)$ and $(N+2\ell-2)^2+4N-4=N^2+4c_\ell$, we have $2C-2=\sqrt{N^2+4\gamma}$ for every block. Hence the order of the levels is the order of the numbers $\gamma$. The smallest value of $\gamma$ is $\gamma_*=2$ for $N\ge3$ and $\gamma_*=1$ for $N=2$, so that the level of a block is $\sqrt{N^2+4\gamma}-\sqrt{N^2+4\gamma_*}$. It is smaller than $4$ if and only if
\begin{equation}\label{threshold}
\gamma<\gamma_*+4+2\sqrt{N^2+4\gamma_*} .
\end{equation}

\emph{The case $N\ge4$.} For the poloidal blocks of degree $k=1,2,3,4$ we have $\gamma=2$, $N+3$, $2N+6$, $3N+11$, and for the toroidal blocks of degree $\ell=1,2,3,4$ we have $\gamma=N-1$, $2N$, $3N+3$, $4N+8$. The numbers $\gamma$ increase with the degree. Since $N\ge4$,
\[
2<N-1<N+3<2N<2N+6<3N+3<3N+11<4N+8 .
\]
The right-hand side of \eqref{threshold} is $6+2\sqrt{N^2+8}$. We have $2N+6<6+2\sqrt{N^2+8}$. Next, $3N+3<6+2\sqrt{N^2+8}$ if and only if $(3N-3)^2<4(N^2+8)$, that is, $5N^2-18N-23<0$. This holds for $N=4$, since $80-72-23<0$, and it fails for $N\ge5$, since $5N^2-18N-23$ is increasing for $N\ge2$ and equals $12$ at $N=5$. Finally, $3N+11>6+2\sqrt{N^2+8}$, since $(3N+5)^2>4(N^2+8)$. This gives (i) and (ii). For $N=4$ we use $\sqrt{28}=2\sqrt7$, $\sqrt{24}=2\sqrt6$, $\sqrt{44}=2\sqrt{11}$, $\sqrt{48}=4\sqrt3$, $\sqrt{72}=6\sqrt2$ and $\sqrt{76}=2\sqrt{19}$ in Theorem \ref{thm:sum}(ii).

\emph{The case $N=3$.} Here both blocks of degree $k$ have $\gamma=k(k+1)$, which takes the values $2,6,12,20,\dots$. The right-hand side of \eqref{threshold} is $6+2\sqrt{17}$, which lies between $14$ and $15$, since $16<17<\frac{81}{4}$. This gives (iii).

\emph{The case $N=2$.} Here $\gamma=1$ for $\bU_T^{(1)}$, and $\gamma=k^2+1$ for $\bU_{P,k}$, which takes the values $2,5,10,17,\dots$. The right-hand side of \eqref{threshold} is $5+4\sqrt2$, which lies between $10$ and $11$, since $25<32<36$. This gives (iv).
\end{proof}

\subsection{Proof of Theorem \ref{thm:sum}}
By Lemma \ref{lem:levels}, each estimate of Theorem \ref{thm:sum} has the form
\begin{equation}\label{sum}
\delta_{2,N}(\bU)\ \ge\ \sum_{i=1}^{n}(\tau_i-\tau_{i-1})\dist^2(\bU,\cG_{i-1}) ,
\end{equation}
where $0=\tau_0<\tau_1<\dots<\tau_{n-1}$ are the levels of the blocks smaller than $4$, $\tau_n=4$, and $\cG_i$ is the sum of the spaces $\cE(\bV)$ of the blocks $\bV$ with level at most $\tau_i$. For instance, for $N\ge5$ we have $\cG_0=\cE_1$, $\cG_1=\cE_1+\cT_1$, $\cG_2=\cE_2+\cT_1$, $\cG_3=\cE_2+\cT_2$ and $\cG_4=\cE_3+\cT_2$. By \eqref{distsplit}, the right-hand side of \eqref{sum} is a sum over the blocks. We compute the contribution of each block $\bV$ and compare it with $\delta_{2,N}(\bV)$.

Let $\bV$ be a block with level $\tau_j$ for some $j\le n-1$. By \eqref{distsplit}, $\bV$ contributes $\|\bV\|^2$ to $\dist^2(\bU,\cG_{i-1})$ for $i\le j$, and $\dist^2(\bV,\cE(\bV))$ for $i>j$. Its contribution to the right-hand side of \eqref{sum} is therefore
\[
\tau_j\|\bV\|^2+(4-\tau_j)\dist^2(\bV,\cE(\bV))\ \le\ \tau_j\|\bV\|^2+4\dist^2(\bV,\cE(\bV))\ \le\ \delta_{2,N}(\bV)
\]
by \eqref{blockest}. Let $\bV$ be a block with level larger than $4$. It contributes $\|\bV\|^2$ to every distance, and its contribution is $4\|\bV\|^2$, which is at most $\delta_{2,N}(\bV)$ by \eqref{blockest}. Summing over the blocks, we obtain \eqref{sum}.

If $\bU\in\cG_{n-1}$, then $\bV\in\cE(\bV)$ for every block with level smaller than $4$, and the other blocks vanish. By the equality cases of Proposition \ref{prop:pol} and \eqref{torall}, $\delta_{2,N}(\bV)=\tau_j\|\bV\|^2$ for a block with level $\tau_j$, which is its contribution. Hence equality holds in \eqref{sum}.

For the sharpness, we show that for every $i\in\{1,\dots,n\}$ there is a solenoidal field $\bE\ne0$ which gives equality in \eqref{sum}, such that the $i$-th term of the right-hand side is nonzero and the terms of index larger than $i$ vanish. Then the $i$-th constant cannot be increased while the other ones are kept fixed. Let first $i\le n-1$. We take a nonzero field $\bE$ in the space $\cE(\bV)$ of a block with level $\tau_i$. Such a field exists, by the end of Section \ref{sec:tor} for the toroidal blocks. Then $\bE\in\cG_i$, so that equality holds, and $\dist(\bE,\cG_{i-1})=\|\bE\|$ since $\bE$ is orthogonal to $\cG_{i-1}$. Let now $i=n$. We take for $\bE$ a radial excitation of an extremal, that is, the field $\bR$ of Section \ref{sec:pol} with $k=1$ if $N\ge3$, and the field $\bR$ of Section \ref{sec:tor} if $N=2$. By \eqref{Rpol} and \eqref{Rtor}, $\delta_{2,N}(\bE)=4\|\bE\|^2$, since $C_N=C_P(N,1)$ for $N\ge3$ and $C_2=C_T(2,1)$. Moreover, $\bE$ is orthogonal in $\dot H^1$ to $\cG_{n-1}$. Indeed, $\bE$ lies in the block $\bU_{P,1}$, respectively $\bU_T^{(1)}$, it is orthogonal to the space $\cE$ of this block, and the other spaces lie in the other blocks. Hence all the distances are equal to $\|\bE\|$, and the right-hand side of \eqref{sum} is $\tau_n\|\bE\|^2=4\|\bE\|^2=\delta_{2,N}(\bE)$. This proves Theorem \ref{thm:sum}. \qed

\section{The stability with the deficit \texorpdfstring{$\delta_{1,N}$}{delta(1,N)}--Proof of Theorem \texorpdfstring{\ref{thm:prod}}{1.3}}\label{sec:prod}

We use the dilation argument of \cite[Lemma 5.1]{DDLLNsol}. Let $\bU\ne0$ be a solenoidal field. Then $\int|\Delta\bU|^2$ and $\int|\nabla\bU|^2|x|^2$ are positive, since if one of them vanishes, then $\nabla\bU=0$. By \eqref{scaling},
\[
\int_{\R^N}|\Delta\bU_{\lambda_*}|^2dx=\int_{\R^N}|\nabla\bU_{\lambda_*}|^2|x|^2dx=\Big(\int_{\R^N}|\Delta\bU|^2dx\int_{\R^N}|\nabla\bU|^2|x|^2dx\Big)^{\frac12}.
\]
Since $\|\bU_{\lambda_*}\|=\|\bU\|$, this gives
\begin{equation}\label{transfer}
\delta_{1,N}(\bU)=\frac12\,\delta_{2,N}(\bU_{\lambda_*}).
\end{equation}
Moreover, the dilations are isometries for $\|\cdot\|$, so that $\dist(\bU_{\lambda_*},\cE)=\dist(\bU,\{\bE_{1/\lambda_*}:\bE\in\cE\})\ge\dist(\bU,\cE^{\rm dil})$ for every family $\cE$.

We also need to know that the test fields of the proof of Theorem \ref{thm:sum} are balanced.

\begin{lemma}\label{lem:balanced}
Let $\bE$ be a nonzero field of one of the spaces $\cE(\bV)$, or one of the radial excitations $\bR$ of Sections \ref{sec:pol} and \ref{sec:tor}. Then
\[
\int_{\R^N}|\Delta\bE|^2dx=\int_{\R^N}|\nabla\bE|^2|x|^2dx .
\]
\end{lemma}

\begin{proof}
We write $I_1=\int|\Delta\bE|^2$ and $I_2=\int|\nabla\bE|^2|x|^2$. A field of $\cE(\bU_{P,k})$ is a sum of fields $\bD(r^{k-1}p_k\phi_{k,j})$ with the same profile, and $I_1$, $I_2$ split over $j$ by the orthogonality (c) of Subsection \ref{sec:sol}. So it is enough to consider one spherical harmonic. Let first $\bE=\bD(f_k\phi_k)$ be a poloidal field whose function $F_k$ is $v(|\xi|)e^{-|\xi|^2/2}$, with $v(\rho)=q(\rho^2)$ and $q$ a polynomial of degree at most one, and let $C=C_P(N,k)$. By Proposition \ref{prop:lift} and Plancherel's identity, $I_1$ is a positive multiple of $\int|\cF G_k|^2|\xi|^6d\xi=\int|F_k|^2|\xi|^{6+2\beta_k}d\xi$. By \eqref{blockw}, \eqref{normw} and \eqref{gauss}, there is a constant $c>0$ such that
\begin{gather}
\|\bE\|^2=c\int_0^\infty|v|^2\rho^{2C-1}e^{-\rho^2}d\rho,\label{bal}\\
I_1=c\int_0^\infty|v|^2\rho^{2C+1}e^{-\rho^2}d\rho,\notag\\
I_1+I_2-2C\|\bE\|^2=c\int_0^\infty|v'|^2\rho^{2C-1}e^{-\rho^2}d\rho .\notag
\end{gather}
Let next $\bE$ be a toroidal field of degree $\ell$ whose functions $w_{ij}$ in \eqref{torblockW} are $a_{ij}\rho^{D_\ell/2}v(\rho)$, with constants $a_{ij}$, $v(\rho)=q(\rho^2)$ and $q$ a polynomial of degree at most one, and let $C=C_T(N,\ell)$. Since $I_1$ is $(2\pi)^{-N}$ times the integral of $|\bW_T|^2|\xi|^4e^{-|\xi|^2}$, where $\bW_T$ is defined from $\bE$ as in \eqref{torw}, the identities \eqref{torw}, \eqref{normv} and \eqref{subst} show that the identities \eqref{bal} hold again, for some constant $c>0$. We pass to the variable $t=\rho^2$, so that $\rho^{2C-1}d\rho=\frac12t^{C-1}dt$ and $|v'|^2=4t|q'|^2$. If $q$ is constant, we may take $q=1$, and then $\|\bE\|^2=\frac c2\Gamma(C)$ and $I_1=\frac c2\Gamma(C+1)=C\|\bE\|^2$, so that $I_1+I_2=2C\|\bE\|^2$ and $I_2=I_1$. If $q(t)=C-t$, a direct computation gives
\[
\int_0^\infty(C-t)^2t^{C-1}e^{-t}dt=C\,\Gamma(C),\qquad \int_0^\infty(C-t)^2t^{C}e^{-t}dt=C(C+2)\Gamma(C).
\]
Hence $I_1=(C+2)\|\bE\|^2$. Moreover $2c\int_0^\infty t^{C}e^{-t}dt=2c\,C\,\Gamma(C)=4\|\bE\|^2$, so that $I_1+I_2=(2C+4)\|\bE\|^2$ and again $I_2=I_1$. The fields of the spaces $\cE(\bV)$ have constant $q$, and the radial excitations have $q$ proportional to $C-t$, in both the poloidal and the toroidal case. This proves the lemma.
\end{proof}

\begin{proof}[Proof of Theorem \ref{thm:prod}]
By \eqref{transfer} and Theorem \ref{thm:sum} applied to $\bU_{\lambda_*}$, we get Theorem \ref{thm:prod}(i). Theorem \ref{thm:prod}(ii) follows, since $\dist(\bU_{\lambda_*},\cE)\ge\dist(\bU,\cE^{\rm dil})$ for every family $\cE$.

For the sharpness in Theorem \ref{thm:prod}(i), we use the test fields $\bE$ of the proof of Theorem \ref{thm:sum}. By Lemma \ref{lem:balanced}, $\lambda_*(\bE)=1$. Hence $\bE_{\lambda_*}=\bE$ and $\delta_{1,N}(\bE)=\frac12\delta_{2,N}(\bE)$, so that the equality in \eqref{sum} for $\bE$ becomes an equality in Theorem \ref{thm:prod}(i), with the same nonzero terms. Hence none of the constants can be increased while the other ones are kept fixed.

For Theorem \ref{thm:prod}(ii), let $\bE$ be a test field of the proof of Theorem \ref{thm:sum} which lies in the space $\cE(\bV)$ of a block with level $\tau_i<4$. The dilations preserve the poloidal fields of each degree and the toroidal fields of each degree. Hence the dilation orbits of the spaces of the blocks with smaller levels are orthogonal to $\bE$ in $\dot H^1$, and $\dist(\bE,\cG_j^{\rm dil})=\|\bE\|$ for $j<i$. On the other hand, $\dist(\bE,\cG_j^{\rm dil})=0$ for $j\ge i$, since $\bE\in\cG_i$. Since $\delta_{1,N}(\bE)=\frac12\delta_{2,N}(\bE)=\frac12\tau_i\|\bE\|^2$, equality holds in Theorem \ref{thm:prod}(ii). Hence all the constants except the last one cannot be increased while the other ones are kept fixed.
\end{proof}

It remains to characterize the equality cases in Theorem \ref{thm:sharp}. Let $\cG_0$ be the first family in Theorem \ref{thm:sum}, that is, $\cE_1$ for $N\ge4$, $\cE_1+\cT_1$ for $N=3$, and $\cT_1$ for $N=2$. If $\bE\in\cG_0$, then $\delta_{2,N}(\bE)=0$ by Theorem \ref{thm:sum}, and hence $0\le2\delta_{1,N}(\bE)\le\delta_{2,N}(\bE)=0$. By the dilation invariance of $\delta_{1,N}$, equality holds in \eqref{2ndHUPsol} on $\cG_0^{\rm dil}$. Conversely, if $\bU\ne0$ and $\delta_{1,N}(\bU)=0$, then $\dist(\bU_{\lambda_*},\cG_0)=0$ by Theorem \ref{thm:prod}(i), since its first constant is positive. Since $\cG_0$ is a finite dimensional subspace, $\bU_{\lambda_*}\in\cG_0$, and $\bU\in\cG_0^{\rm dil}$. This completes the proof of Theorem \ref{thm:sharp}. \qed

\medskip
\noindent\textbf{Acknowledgements.} A. T. Duong and V. H. Nguyen were supported by the Vietnam National Foundation for Science and Technology Development (NAFOSTED) under grant number 101.02-2025.33. N. Lam was partially supported by an NSERC Discovery Grant. G. Lu was partially supported by grants from the Simons Foundation.

\medskip

{\bf AI Declaration.} The mathematical content of this paper, including all results and proofs, is the work of the authors, who take full responsibility for the content of the paper. AI tools were used only to improve the language of the manuscript after the work was completed.


\end{document}